\documentclass[final, nohypdvips]{siamart0516}
\pdfoutput=1

\usepackage{braket,amsfonts}

\usepackage{array}

\usepackage[caption=false]{subfig}

\usepackage{graphicx}
\usepackage{pst-pdf}
\usepackage{comment}
\usepackage{multirow}

\newsiamthm{claim}{Claim}
\newsiamremark{rem}{Remark}
\newsiamremark{expl}{Example}
\newsiamremark{hypothesis}{Hypothesis}
\crefname{hypothesis}{Hypothesis}{Hypotheses}

\usepackage{algorithmic}

\usepackage{graphicx,epstopdf}

\newcommand{\RR}[1]{\mathbb{R}^{#1}}
\usepackage{mathtools}
\newcommand{\defeq}{\coloneqq}
\newcommand{\Span}[1]{\text{span}\{#1\}}
\newcommand{\SG}{^\text{SG}}

\newcommand{\PS}{^\text{PS}}
\newcommand{\LSPG}{^\text{LSPG}}
\newcommand{\wLSPG}{^{\text{LSPG}(M)}}  
\newcommand{\wLSPGlarger}{^{\text{LSPG}'(M)}} 
\newcommand{\wLSPGArg}[1]{^{\text{LSPG}(#1)}}
\newcommand{\RHS}[1]{b#1}
\newcommand{\numxi}{{n_\xi}}
\newcommand{\numpsi}{{n_\psi}}
\newcommand{\LSV}[1]{\sigma_{\max}(#1)} 
\newcommand{\SSV}[1]{\sigma_{\min}(#1)} 

\newcommand{\LSPGSG}{LSPG($A$)/SG}  
\newcommand{\LSPGPS}{LSPG($2$)/PS}
\newcommand{\LSPGI}{LSPG($A^TA$)}
\newcommand{\LSPGQ}{LSPG($F^TF$)}

\newcommand{\SGsup}{\wLSPGArg{A}} 
\newcommand{\PSsup}{\wLSPGArg{2}}
\newcommand{\LSPGsup}{\wLSPGArg{A^TA}}

\newcommand{\errorsol}[1]{{e}(#1)}

\usepackage{dsfont}
\usepackage{amsmath}
\usepackage{amssymb}

\Crefname{ALC@unique}{Line}{Lines}

\usepackage{amsopn}

\title{Stochastic Least-squares Petrov--Galerkin Method for Parameterized Linear Systems%
  \thanks{This work was supported by the U.S. Department of Energy Office of Advanced Scientific Computing Research, Applied Mathematics program under award DEC-SC0009301 and by the U.S. National Science Foundation under grant DMS1418754.}}

\author{Kookjin Lee%
  \thanks{Department of Computer Science, University of Maryland, College Park, MD 20742 (\email{klee@cs.umd.edu}).}%
  \and
  Kevin Carlberg%
  \thanks{Sandia National Laboratories, Livermore, CA 94550
    (\email{ktcarlb@sandia.gov}).
	Sandia National Laboratories is a multi-program laboratory managed and
	operated by Sandia Corporation, a wholly owned subsidiary of Lockheed Martin
	Corporation, for the U.S. Department of Energy's National Nuclear Security
	Administration under contract DE-AC04-94AL85000.}
  \and
  Howard C. Elman%
  \thanks{Department of Computer Science and Institute for Advanced Computer Studies, University of Maryland, College Park, MD 20742
    (\email{elman@cs.umd.edu}).}
}

\headers{Stochastic Least-Squares Petrov--Galerkin Method}
{K.\ Lee, K.\ Carlberg, and H.\ C.\ Elman}

\begin{document}
\maketitle

\begin{abstract}
We consider the numerical solution of parameterized linear
systems where the system matrix, the solution, and the right-hand side are
parameterized by a set of uncertain input parameters. We explore spectral
methods in which the solutions are approximated in a chosen finite-dimensional subspace. It has been shown that the stochastic Galerkin projection technique fails to minimize any measure of the solution error \cite{mugler2013convergence}.
As a remedy for this, we propose a novel stochastic least-squares Petrov--Galerkin (LSPG) method. The proposed method is
optimal in the sense that it produces the solution that minimizes a weighted
$\ell^2$-norm of the residual over all solutions in a given
finite-dimensional subspace. Moreover, the method
can be adapted to minimize the solution error in different weighted
$\ell^2$-norms by simply applying a weighting function within the least-squares
formulation. In addition, a goal-oriented semi-norm induced by an
output quantity of interest can be minimized by defining a weighting function
as a linear functional of the solution. We establish optimality and error
bounds for the proposed method, and extensive numerical experiments show that
the weighted LSPG methods outperforms other spectral methods in
minimizing  corresponding target weighted norms.
\end{abstract}

\begin{keywords}
stochastic Galerkin, least-squares Petrov--Galerkin projection, residual minimization, spectral projection
\end{keywords}

\begin{AMS}
35R60, 60H15, 60H35, 65C20, 65N30, 93E24
\end{AMS}

\section{Introduction} \label{sec:intro}
Forward uncertainty propagation for parameterized linear systems is important
in a range of applications to characterize the effects of uncertainties on the
output of computational models. Such parameterized linear systems arise in
many important problems in science and engineering, including
stochastic partial differential equations (SPDEs) where uncertain input parameters are modeled as a set of
random variables (e.g., diffusion/ground water flow simulations where
diffusivity/permeability is modeled as a random field \cite{holden1996stochastic, zhang2001stochastic}). It has
been shown \cite{elman2011assessment} that intrusive methods (e.g., stochastic Galerkin
\cite{babuska2004galerkin, deb2001solution, ghanem2003stochastic,
matthies2005galerkin, xiu2002modeling}) for uncertainty propagation can lead
to smaller errors for a fixed basis dimension, compared with non-intrusive
methods \cite{xiu2010numerical} (e.g., sampling-based methods \cite{barth2011multi,
graham2011quasi, kuo2012quasi}, stochastic collocation
\cite{babuvska2007stochastic, nobile2008sparse}). 

The stochastic Galerkin method combined with generalized polynomial chaos
(gPC) expansions  \cite{xiu2002wiener} seeks a polynomial approximation of the
numerical solution in the stochastic domain by enforcing a Galerkin
orthogonality condition, i.e., the residual of the parameterized linear system
is forced to be orthogonal to the span of the stochastic polynomial basis with
respect to an inner product associated with an underlying probability measure.
The Galerkin projection scheme is popular for its simplicity (i.e., the trial
and test bases are the same) and its optimality in terms of minimizing an energy norm of solution
errors when the underlying PDE operator is symmetric positive definite.  In many applications, however, it has been shown that the stochastic Galerkin method does not exhibit any optimality property \cite{mugler2013convergence}. That is, it does not produce solutions that minimize any measure of the solution error. In such cases, the stochastic Galerkin method can lead to poor approximations and non-convergent behavior.

To address this issue, we propose a novel optimal projection technique, which we refer to as the stochastic least-squares Petrov--Galerkin (LSPG) method. Inspired by the successes of LSPG methods in nonlinear model reduction \cite{CarlbergGappy, carlberg2013gnat, carlbergGalDiscOpt}, finite element methods \cite{bochev1998finite, bochev2009least,jiang1990least}, and iterative linear solvers (e.g., GMRES, GCR) \cite{saadBook}, we propose, as an alternative to enforcing the Galerkin orthogonality condition, to directly minimize the residual of a parameterized linear system over the stochastic domain in a (weighted) $\ell^2$-norm. The stochastic LSPG method produces an optimal solution for a given stochastic subspace and guarantees that the $\ell^2$-norm of the residual monotonically decreases as the stochastic basis is enriched. In addition to producing monotonically convergent approximations as measured in the chosen weighted $\ell^2$-norm, the method can also be adapted to target output quantities of interest (QoI); this can be accomplished by employing a weighted $\ell^2$-norm used for least-squares minimization that coincides with the $\ell^2$-(semi)norm of the error in the chosen QoI. 

In addition to proposing the stochastic LSPG method, this study shows that specific choices of weighting functions lead to equivalence between the stochastic LSPG method and both the stochastic Galerkin
method and the pseudo-spectral method  \cite{xiu2007efficient, xiu2010numerical}. 
We demonstrate the effectiveness of this method with extensive numerical
experiments on various SPDEs. The results show that the proposed LSPG
technique significantly outperforms the stochastic Galerkin when the solution
error is measured in different weighted $\ell^2$-norms. We also show that the
proposed method can effectively minimize the error in target QoIs.

An outline of the paper is as follows.  Section \ref{sec:spectral_methods} formulates parameterized linear algebraic systems and reviews conventional spectral approaches for computing numerical solutions. Section \ref{sec:slspg} develops a residual minimization formulation based on least-squares methods and its adaptation to the stochastic LSPG method. We also provide proofs of optimality and monotonic convergence behavior of the proposed method. Section \ref{sec:analysis} provides error analysis for stochastic LSPG methods. Section \ref{sec:results} demonstrates the efficiency and the effectiveness of the proposed methods by testing them on various benchmark problems. Finally, Section \ref{sec:conclusion} outlines some conclusions.

\section{Spectral methods for parameterized linear systems}\label{sec:spectral_methods}
We begin by introducing a mathematical formulation of parameterized linear
systems and briefly reviewing the stochastic Galerkin and the
{pseudo-spectral} methods 
, which are spectral methods for approximating the numerical solutions of such systems. 

\subsection{Problem formulation}\label{sec:problem_form}
Consider a parameterized linear system
\begin{equation}\label{eq:param_sys}
A(\xi) u(\xi) = \RHS{(\xi)},
\end{equation}
where $A:\Gamma\rightarrow\RR{n_x\times n_x}$, and $u,b: \Gamma\rightarrow\RR{n_x}$.  
The system is parameterized by a set of stochastic input parameters
$\xi(\omega) \equiv \{ \xi_1(\omega),\ldots,\xi_\numxi(\omega)\}$. Here, $\omega\in\Omega$
is an elementary event in a probability space  $(\Omega,\mathcal F, P)$ and
the stochastic domain is denoted by $\Gamma \equiv \prod_{i=1}^{\numxi}
\Gamma_i$ where $\xi_i: \Omega \rightarrow \Gamma_i$. We are interested in
computing a spectral approximation of the numerical solution $u(\xi)$ in an
$\numpsi$-dimensional
subspace $S_{\numpsi}$ spanned by a finite set of polynomials
$\{\psi_i(\xi)\}_{i=1}^{\numpsi}$ such that
$S_{\numpsi} \equiv \Span{\psi_i}_{i=1}^{\numpsi}\subseteq L^2(\Gamma)$, i.e.,
\begin{equation}\label{eq:spectral_approx}
u(\xi) \approx \tilde{u}(\xi) = \sum_{i=1}^{\numpsi} \bar{u}_i \psi_i(\xi)=(\psi^T(\xi)\otimes I_{n_x})\bar u,
\end{equation}
where  $\{\bar{u}_i\}_{i=1}^{\numpsi}$
with $\bar{u}_i \in \mathbb{R}^{n_x}$ are unknown coefficient vectors,
$\bar{u} \equiv  [\bar u_1^T\ \cdots\ \bar u_{\numpsi}^T]^T \in \mathbb{R}^{n_x \numpsi}$ is the vertical concatenation of these coefficient vectors, $\psi \equiv [\psi_1\ \cdots\ \psi_{\numpsi}]^T\in\RR{\numpsi}$ is a concatenation of the polynomial basis, $\otimes$ denotes the Kronecker product, and $I_{n_x}$ denotes the identity matrix of dimension $n_x$. Note
that $\tilde{u}\in(S_{\numpsi})^{n_x}$. Typically, the ``stochastic'' basis $\{\psi_i\}$ consists of products of univariate polynomials: $\psi_i \equiv \psi_{\alpha(i)} \equiv \prod_{k=1}^{\numxi} \pi_{\alpha_k(i)}(\xi_k)$ where $\{\pi_{\alpha_k(i)}\}_{k=1}^{\numxi} $ are univariate polynomials, ${\alpha}(i) = (\alpha_1(i),\cdots,\alpha_{\numxi}(i)) \in \mathbb{N}^{\numxi}_0$ is a multi-index and $\alpha_k$ represents the degree of a polynomial in $\xi_k$. The dimension of the stochastic subspace
$\numpsi$ depends on the number of random variables $\numxi$, the maximum polynomial
degree  $p$, and a construction of the polynomial space (e.g., a total-degree
space that contains polynomials with total degree up to $p$, $\sum_{k=1}^{\numxi}\alpha_k(i) \leq p$). By substituting
$u(\xi)$ with $\tilde{u}(\xi)$ in \eqref{eq:param_sys}, the residual can be defined
as
\begin{align}\label{eq:residual_def}
r(\bar{u};\xi) &\defeq \RHS{(\xi)}- A(\xi) \sum_{i=1}^{\numpsi} \bar{u}_i \psi_i(\xi)=\RHS{(\xi)} -  (\psi^T(\xi) \otimes A(\xi)) \bar{u},
\end{align}
where $\psi^T(\cdot) \otimes A(\cdot): \Gamma \rightarrow \mathbb{R}^{n_x \times \numpsi n_x}$. 

It follows from \eqref{eq:spectral_approx} and \eqref{eq:residual_def} that our goal now is to compute the unknown coefficients $\{\bar{u}_i\}_{i=1}^{\numpsi}$ of the solution expansion.
We briefly review two conventional approaches for doing so: the stochastic Galerkin method and the pseudo-spectral method. In the following, $\rho \equiv \rho(\xi)$ denotes an underlying measure of the stochastic space $\Gamma$ and
\begin{align}
\langle g, h \rangle_\rho &\equiv \int_\Gamma g(\xi) h(\xi) \rho(\xi) d\xi, \label{eq:def_inner_prod}\\
E[ g ] &\equiv \int_\Gamma g(\xi) \rho(\xi) d\xi, \label{eq:def_expectation}
\end{align}
define an inner product between scalar-valued functions $g(\xi)$ and $h(\xi)$ with respect to $\rho(\xi)$ and the expectation of $g(\xi)$, respectively. The $\ell^2$-norm of a vector-valued function $v(\xi) \in \RR{n_x}$ is defined as
\begin{align}\label{eq:def_l2_norm}
\| v \|_2^2 & \equiv \sum_{i=1} ^{n_x} \int_\Gamma v^2_i(\xi) \rho(\xi)d\xi = E[v^T v] .
\end{align}
Typically, the polynomial basis is constructed to be orthogonal in the $\langle\cdot,\cdot\rangle_\rho$ inner product, i.e., $\langle\psi_i,\psi_j\rangle_\rho = \prod_{k=1}^{\numxi} \langle \pi_{\alpha_k(i)}, \pi_{\alpha_k(j)} \rangle_{\rho_k} = \delta_{ij}$,
where $\delta_{ij}$ denotes the Kronecker delta. 

\subsection{Stochastic Galerkin method} \label{sec:SG_method}
The stochastic Galerkin method computes the unknown coefficients $\{\bar{u}_i\}_{i=1}^{\numpsi}$ of $\tilde{u}(\xi)$ in \eqref{eq:spectral_approx} by imposing orthogonality of the residual \eqref{eq:residual_def} with respect to the inner product $\langle \cdot,\cdot\rangle_\rho$ in the subspace $S_{\numpsi}$.  This Galerkin orthogonality condition can be expressed as follows: Find $\bar{u}\SG\in\RR{n_x\numpsi}$ such that
\begin{equation}\label{eq:galerkin_condition}
\langle r_i(\bar{u}\SG), \psi_j \rangle_{\rho} = E[ r_i(\bar{u}\SG) \psi_j ] = 0, \qquad i=1,\ldots,n_x,\ j = 1,\ldots,\numpsi,
\end{equation}
where $r\equiv [r_1\ \cdots\ r_{n_x}]^T$. The condition \eqref{eq:galerkin_condition} can be represented in matrix notation as
\begin{equation}\label{eq:galerkin_matrix_form}
 E[ \psi \otimes r(\bar{u}\SG) ] = 0.
\end{equation}
From the definition of the residual \eqref{eq:residual_def}, this gives a system of linear equations
\begin{equation}\label{eq:sg_linear_sys}
E[\psi \psi^T \otimes A] \bar{u}\SG = E[\psi \otimes \RHS{}],
\end{equation}
of dimension $n_x \numpsi$. This yields an algebraic expression for the stochastic-Galerkin approximation
\begin{equation}\label{eq:sg_linear_sys_proj}
\tilde u\SG(\xi) = (\psi(\xi)^T\otimes I_{n_x})E[\psi\psi^T\otimes A]^{-1}E[\psi\otimes A u].
\end{equation}
If $A(\xi)$ is symmetric positive definite, the solution
of linear system \eqref{eq:sg_linear_sys} minimizes the solution error
$\errorsol{x} \equiv u - x$ in the $A(\xi)$-induced energy norm $\| v  \|_A^2 \equiv  E[ v^T A v]$, i.e.,
\begin{equation}\label{eq:stochGalOpt}
\tilde{u}\SG(\xi) = \underset{x\in (S_{\numpsi})^{n_x}}{\arg\min}\|\errorsol{x}\|_A^2.
\end{equation}
In general, however, the stochastic-Galerkin approximation does not minimize any measure
of the solution error.  


\subsection{Pseudo-spectral method} \label{sec:PS_method}
The pseudo-spectral method directly approximates the unknown coefficients $\{\bar{u}_i\}_{i=1}^{\numpsi}$ of $\tilde{u}(\xi)$ in \eqref{eq:spectral_approx} by exploiting orthogonality of the polynomial basis $\{ \psi_i(\xi)\}_{i=1}^{\numpsi}$. That is, the coefficients $\bar{u}_i$ can be obtained by projecting the numerical solution $u(\xi)$ onto the orthogonal polynomial basis as
\begin{equation}\label{eq:ps_projection}
\bar{u}_i\PS = E[ u \psi_i ],   \quad i=1,\ldots,\numpsi, 
\end{equation}
which can be expressed as
\begin{equation}
\bar{u}\PS = E[ \psi \otimes A^{-1} \RHS{}],
\end{equation}
or equivalently
\begin{equation}\label{eq:ps_linear_sys_proj}
\tilde u\PS(\xi) = (\psi(\xi)^T\otimes I_{n_x})E[\psi\otimes u].
\end{equation}
The associated optimality property of the approximation, which can be derived from optimality of orthogonal projection, is
\begin{equation}\label{eq:spectralProjOpt}
\tilde{u}\PS(\xi) = \underset{x\in (S_{\numpsi})^{n_x}}{\arg\min}\|\errorsol{x}\|_2^2.
\end{equation}
In practice, the coefficients $\{\bar{u}\PS_i\}_{i=1}^{\numpsi}$ are approximated via numerical quadrature as
\begin{align}\label{eq:ps_proj_quad}
\bar{u}_i\PS = E[ u \psi_i ]=\sum_{k=1}^{n_q} u(\xi^{(k)}) \psi_i(\xi^{(k)}) w_k= \sum_{k=1}^{n_q} \left( A^{-1}(\xi^{(k)}) f(\xi^{(k)}) \right) \psi_i(\xi^{(k)}) w_k,
\end{align}
where  $\{(\xi^{(k)}, w_k)\}_{k=1}^{n_q}$ are the quadrature points and weights.

While stochastic Galerkin leads to an optimal approximation \eqref{eq:stochGalOpt} under certain conditions and pseudo-spectral projection minimizes the $\ell^2$-norm of the solution error \eqref{eq:spectralProjOpt}, neither approach provides the flexibility to tailor the optimality properties of the approximation. This may be important in applications where, for example, minimizing the error in a quantity of interest is desired. To address this, we propose a general optimization-based framework for spectral methods that enables the choice of a targeted weighted $\ell^2$-norm in which the solution error is minimized. 

\section{Stochastic least-squares Petrov--Galerkin method}\label{sec:slspg}
As a starting point, we propose a residual-minimizing formulation that computes the coefficients $\bar{u}$ by directly minimizing the $\ell^2$-norm of the residual, i.e.,
\begin{align}\label{eq:min_res_L2_first}
\begin{split}
\tilde u\LSPG(\xi) &= \underset{x \in (S_{\numpsi})^{n_x}}{\arg\min}{ \| f - Ax\|_2^2 }
=
\underset{x \in (S_{\numpsi})^{n_x}}{\arg\min}{ \| \errorsol{x}\|_{A^TA}^2 },
\end{split}
\end{align}
where $\|v\|_{A^TA}^2\equiv E[v^TA^TAv]$. Thus, the $\ell^2$-norm of the residual is equivalent to a weighted $\ell^2$-norm of the solution error.  
Using \eqref{eq:spectral_approx} and
\eqref{eq:residual_def}, we have
\begin{equation}\label{eq:min_res_L2}
\bar{u}\LSPG = \underset{\bar{x} \in \RR{n_x \numpsi}}{\arg\min}{ \| r(\bar{x})\|_2^2 }.
\end{equation}
The definition of the residual \eqref{eq:residual_def} allows the objective function in \eqref{eq:min_res_L2} to be written in quadratic form as
\begin{align}\label{eq:res_quadratic}
\begin{split}
\| r( \bar{x}) \|_2^2  &= \| f- (\psi^T \otimes  A) \bar{x} \|_2^2 = \bar{x}^T E[\psi\psi^T \otimes A^T A ] \bar{x} - 2 E[ \psi \otimes A^T f]^T \bar{x} + E[f^T f].
\end{split}
\end{align}
Noting that the mapping $\bar x\mapsto \| r( \bar{x}) \|_2^2$ is convex, the (unique) solution $\bar{u}\LSPG$ to \eqref{eq:min_res_L2} is a stationary point of $\| r( \bar{x}) \|_2^2$ and thus satisfies 
\begin{equation}\label{eq:res_NE}
E[\psi \psi^T \otimes A^TA] \bar{u}\LSPG = E[\psi \otimes A^T f],
\end{equation}
which can be interpreted as the normal-equations form of the linear least-squares problem \eqref{eq:min_res_L2}.

Consider a generalization of this idea that minimizes the solution error in a
targeted weighted $\ell^2$-norm by choosing a specific weighting function. 
Let us define a weighting function $M(\xi) \equiv M_\xi(\xi) \otimes M_x(\xi) $, where $M_\xi: \Gamma \rightarrow \mathbb{R}$ and $M_x: \Gamma \rightarrow \mathbb{R}^{n_x \times n_x}$. Then, the stochastic LSPG method can be written as 
\begin{align}\label{eq:min_res_L2_first_weight}
\begin{split}
\tilde u\wLSPG(\xi) &= \underset{x \in (S_{\numpsi})^{n_x}}{\arg\min}{ \| M(\RHS{} - Ax)\|_2^2 }
=
\underset{x \in (S_{\numpsi})^{n_x}}{\arg\min}{ \| \errorsol{x}\|_{A^TM^TMA}^2 },
\end{split}
\end{align}
with $\|v\|_{A^TM^TMA}^2\equiv E[v^TA^TM^TMAv]= E[(M_\xi^T M_\xi \otimes
(M_xAv)^T M_x Av]$. Algebraically, this is equivalent to
\begin{align}\label{eq:stochastic_lspg_algebraic}
\begin{split}
\bar{u}\wLSPG &= \underset{{\bar{x} \in \RR{n_x \numpsi} }}{\arg\min}{\| M r( \bar{x}) \|_2^2}
= \underset{{\bar{x} \in \RR{n_x \numpsi} }}{\arg\min}{ \| (M_\xi \otimes M_x) (1 \otimes \RHS{}- \left( \psi^T \otimes A \right) \bar{x}) \|_{2}^2}\\
&= \underset{{\bar{x} \in \RR{n_x \numpsi} }}{\arg\min}{ \|  M_\xi \otimes (M_x \RHS{}) - \left( (M_\xi \psi^T) \otimes (M_x A) \right) \bar{x} \|_{2}^2}.
\end{split}
\end{align}
We will restrict our attention to the case $M_\xi(\xi) = 1$ and denote $M_x(\xi)$ by $M(\xi)$ for simplicity. Now, the algebraic stochastic LSPG problem \eqref{eq:stochastic_lspg_algebraic} simplifies to
\begin{align}\label{eq:min_wlspg}
\bar{u}\wLSPG &= \underset{{\bar{x} \in \RR{n_x \numpsi} }}{\arg\min}{\| M r( \bar{x}) \|_2^2} = \underset{{\bar{x} \in \RR{n_x \numpsi} }}{\arg\min}{\| M \RHS{} - (\psi^T \otimes M A) \bar{x} \|_2^2}.
\end{align}
The objective function in \eqref{eq:min_wlspg} can be written in quadratic form as 
\begin{align}\label{eq:wlspg_quad_form}
\begin{split}
\| Mr(\bar x)\|_2^2 = &\bar{x}^T E[ (\psi \psi^T \otimes MA^TM^T MA)] \bar{x}
- 2 (E[ \psi \otimes A^TM^TMf ])^T \bar{x} + E[\RHS{}^TM^TM\RHS{}].
\end{split}
\end{align}
As before, because the mapping $\bar x\mapsto \|Mr(\bar x)\|_2^2$ is convex, the unique solution $\bar{u}\wLSPG$ of \eqref{eq:min_wlspg} corresponds to a stationary point of $\|Mr(\bar x)\|_2^2$ and thus satisfies 
\begin{equation}\label{eq:wlspg_NE}
E[ \psi \psi^T \otimes A^TM^T MA] \bar{u}\wLSPG = E[ \psi \otimes A^T M^TMf ],
\end{equation}
which is the normal-equations form of the linear least-squares problem \eqref{eq:min_wlspg}.
This yields the following algebraic expression for the stochastic-LSPG approximation:
\begin{equation}\label{eq:lspg_linear_sys_proj}
\tilde u\wLSPG(\xi) = (\psi(\xi)^T\otimes I_{n_x})E[\psi\psi^T\otimes A^TM^TMA]^{-1}E[\psi\otimes A^TM^TMA u].
\end{equation}

\textbf{Petrov--Galerkin projection}. Another way of interpreting the normal equations \eqref{eq:wlspg_NE} is that the
(weighted) residual $M(\xi)r(\bar u\wLSPG; \xi) $ is
enforced to be orthogonal to the subspace spanned by the optimal test basis
$\{\phi_i\}_{i=1}^{\numpsi}$ with $\phi_i(\xi)\defeq \psi_i(\xi) \otimes
M(\xi)A(\xi)$ and $\Span{\phi_i}_{i=1}^{\numpsi}\subseteq L^2(\Gamma)$. That is, this projection is precisely the (least-squares) {Petrov--Galerkin}
projection,
\begin{equation}\label{eq:lspg_projection}
E[ \phi^T  ( \RHS{} - (\psi^T \otimes MA)\bar{u}\wLSPG )] = 0,
\end{equation}
where $\phi(\xi)\equiv[\phi_1(\xi)\ \cdots\ \phi_{\numpsi}(\xi)]$.

\textbf{Monotonic Convergence}. The stochastic least-squares Petrov-Galerkin is monotonically convergent. That is, as the trial subspace $S_{\numpsi}$ is enriched (by adding polynomials to the basis), the optimal value of the convex objective function $\| M r(\bar{u}\wLSPG)\|_2^2$ monotonically decreases. This is apparent from the LSPG optimization problem \eqref{eq:min_res_L2_first_weight}: Defining
\begin{align}\label{eq:min_res_L2_first_weight_superspace}
\begin{split}
\tilde u\wLSPGlarger(\xi) &= \underset{x \in (S_{\numpsi+1})^{n_x}}{\arg\min}{ \| M(\RHS{} - Ax)\|_2^2 }
,
\end{split}
\end{align}
we have $\| M(\RHS{} - A\tilde u\wLSPGlarger)\|_2^2\leq \| M(\RHS{} - A\tilde u\wLSPG)\|_2^2$ (and $\|u - u\wLSPGlarger\|_{A^TM^TMA}$ $\leq \|u - u\wLSPG\|_{A^TM^TMA}$) if $S_{\numpsi}\subseteq S_{\numpsi+1}$.

\textbf{Weighting strategies}. Different choices of weighting function $M(\xi)$ allow LSPG to minimize different measures of the error. We focus on four particular choices:
\begin{enumerate}
\item $M(\xi) = C^{-1}(\xi)$, where $C(\xi)$ is a Cholesky factor of $A(\xi)$,
i.e., $A(\xi) = C(\xi)C^T(\xi)$. This decomposition exists if and only if $A$
is symmetric positive semidefinite. In this case, LSPG minimizes the energy
norm of the solution error $\|\errorsol{{x}} \|_A^2 \equiv \| C^{-1} r(\bar{x}) \|_2^2$ ($= \| e( (\Psi^T \otimes I_{n_x})\bar{x})\|_A^2$) and is mathematically equivalent to the stochastic Galerkin method described in Section \ref{sec:SG_method}, i.e., $\tilde u\wLSPGArg{C^{-1}}=\tilde u\SG$. This can be seen by comparing \eqref{eq:stochGalOpt} and \eqref{eq:min_res_L2_first_weight} with $M=C^{-1}$, as $A^TM^TMA = A$ in this case. 
\item $M(\xi) = I_{n_x}$, where $I_{n_x}$ is the identity matrix of dimension
$n_x$. In this case, LSPG minimizes the $\ell^2$-norm of the residual $  \| \errorsol{ x}\|_{A^TA}\equiv  \|r(\bar{x}) \|_2^2$.
\item $M(\xi) = A^{-1}(\xi)$. In this case, LSPG minimizes the $\ell^2$-norm of solution error $\| \errorsol{x} \|_2^2  \equiv \| A^{-1} r(\bar x)\|_2^2$. This is mathematically equivalent to the pseudo-spectral method described in Section \ref{sec:PS_method}, i.e., $\tilde u\wLSPGArg{A^{-1}}=\tilde u\PS$, which can be seen by comparing \eqref{eq:spectralProjOpt} and \eqref{eq:min_res_L2_first_weight} with $M=A^{-1}$.
\item $M(\xi) = F(\xi)A^{-1}(\xi)$ where $F: \Gamma \rightarrow \mathbb{R}^{n_{o} \times
n_x}$ is a linear functional of the solution associated with a vector of output quantities of interest. In this case, LSPG minimizes the $\ell^2$-norm of the error in the output quantities of interest $\|F \errorsol{x} \|_2^2 \equiv \| FA^{-1}r(\bar{x}) \|_2^2 $.
\end{enumerate}
We again emphasize that two particular choices of the weighting function $M(\xi)$ lead to equivalence between LSPG and existing spectral-projection methods (stochastic Galerkin and pseudo-spectral projection), i.e.,
\begin{equation}\label{eq:equivalence}
\tilde u\wLSPGArg{C^{-1}}=\tilde u\SG,\quad\quad \tilde u\wLSPGArg{A^{-1}}=\tilde u\PS,
\end{equation}
where the first equality is valid (i.e., the Cholesky decomposition $A(\xi) = C(\xi)C^T(\xi)$ can be computed) if and only if $A$ is symmetric positive semidefinite.
Table \ref{tab:wlspg_method_table} summarizes the target quantities to minimize (i.e., $\| \errorsol{x} \|_\Theta^2 \equiv E[ \errorsol{x}^T \Theta \errorsol{x}])$, the corresponding LSPG weighting functions, and the method names LSPG($\Theta$).
\begin{table}[htbp]
\begin{center}\footnotesize
\renewcommand{\arraystretch}{1.3}
\setlength{\tabcolsep}{4pt}
\caption{Different choices for the LSPG weighting function}\label{tab:wlspg_method_table}
\begin{tabular}{| l | l | l | l |}
\hline
\multicolumn{2}{|c|}{Quantity minimized by LSPG} & \multirow{2}{*}{Weighting function} & \multirow{2}{*}{Method name}\\
\cline{1-2}
Quantity & Expression & & \\
\hline
Energy norm of error & $\| \errorsol{x} \|_A^2$ &  $M(\xi)=C^{-1}(\xi)$ & \LSPGSG\\
$\ell^2$-norm of residual & $\|\errorsol{x}\|_{A^TA}^2$ & $M(\xi) = I_{n_x}$ & \LSPGI \\
$\ell^2$-norm of solution error & $\| \errorsol{x} \|_2^2$ & $M(\xi) = A^{-1}(\xi)$ & \LSPGPS\\
$\ell^2$-norm of error in quantities of interest  & $\|F\errorsol{x}\|_{2}^2$  & $M(\xi) = F(\xi)A^{-1}(\xi)$ & \LSPGQ \\
\hline
\end{tabular}
\end{center}
\end{table}

\section{Error analysis}\label{sec:analysis}
If an approximation satisfies an optimal-projection condition
 \begin{equation} 
\tilde u = \underset{x \in (S_{\numpsi})^{n_x}}{\arg\min}{ \|
\errorsol{x}\|_{\Theta}^2 },
 \end{equation} 
then 
 \begin{equation} 
 \| \errorsol{\tilde u}\|_{\Theta}^2  =  \min_{x \in (S_{\numpsi})^{n_x}}\| \errorsol{x}\|_{\Theta}^2.
 \end{equation} 
Using norm equivalence
\begin{equation}\label{eq:norm_equivalence}
\|x\|_{\Theta'}^2\leq C\|x\|_{\Theta}^2,
\end{equation}
we 
can characterize the solution error $e(\tilde u)$ in any alternative norm $\Theta'$
as
\begin{gather}
\| \errorsol{\tilde u}\|_{\Theta'}^2  \leq  C\min_{x \in (S_{\numpsi})^{n_x}}\| \errorsol{x}\|_{\Theta}^2.
 \end{gather} 
Thus, the error in an alternative norm $\Theta'$ is controlled by the 
optimal objective-function value
$\min_{x \in (S_{\numpsi})^{n_x}}\|\errorsol{x}\|_\Theta^2$ (which can be made small if the 
trial space admits accurate solutions) and the stability constant $C$.

Table \ref{tab:norm_equivalence} reports norm-equivalence constants for the
norms considered in this work. Here,
we have defined
 \begin{equation} 
 \sigma_\text{min}(M) \equiv \inf_{x\in(L^2(\Gamma))^{n_x}}\|Mx\|_2/\|x\|_2,\quad 
\sigma_\text{max}(M) \equiv \sup_{x\in(L^2(\Gamma))^{n_x}}\|Mx\|_2/\|x\|_2.
  \end{equation} 
\begin{table}[ht] 
\begin{center}\footnotesize
\renewcommand{\arraystretch}{1.5}
\setlength{\tabcolsep}{4pt}
\caption{Stability constant $C$ in \eqref{eq:norm_equivalence}}
 \begin{tabular}{|l|c|c|c|c|c|} 
\hline
&$\Theta'=A$&$\Theta'=A^TA$&$\Theta'=2$ & $\Theta'=F^TF$\\
\hline
$\Theta=A$& 1& $\sigma_\text{max}(A)$& $\frac{1}{\sigma_\text{min}(A)}$& $\frac{\sigma_\text{max}(F)^2}{\sigma_\text{min}(A)}$\\
$\Theta=A^TA$& $\frac{1}{\sigma_\text{min}(A)}$& 1& $\frac{1}{\sigma_\text{min}(A)^2}$& $\frac{\sigma_\text{max}(F)^2}{\sigma_\text{min}(A)^2}$\\
$\Theta=2$& $\sigma_\text{max}(A)$& $\sigma_\text{max}(A)^2$& 1& $\sigma_\text{max}(F)^2$\\
$\Theta=F^TF$& $\frac{\sigma_\text{max}(A)}{\sigma_\text{min}(F)^2}$
& $\frac{\sigma_\text{max}(A)^2}{\sigma_\text{min}(F)^2}$ & $\frac{1}{\sigma_\text{min}(F)^2}$& 1\\
  \hline 
  \end{tabular} 
	\label{tab:norm_equivalence}
\end{center}
  \end{table} 
  
\noindent This exposes several interesting conclusions. First, if $n_o
< n_x$, then the null space of $F$ is nontrivial and so $\sigma_\text{min}(F)
= 0$. This implies that \LSPGQ, for which $\Theta=F^TF$, will have
an undefined value of $C$ when the solution error is measured in other norms, i.e.,
for $\Theta' = A$, $\Theta' = A^TA$, and $\Theta' = 2$. It will have
controlled errors only for $\Theta' = F^TF$, in which case $C=1$. 
Second, note that for problems with small
$\sigma_\text{min}(A)$, the $\ell^2$ norm in the quantities of interest may
be large for the \LSPGSG{}, or \LSPGI{}, while it will remain well behaved for \LSPGPS{} and 
\LSPGQ{}.


\section{Numerical experiments}\label{sec:results}
This section explores the performance of the LSPG
methods for solving elliptic SPDEs parameterized by one random variable (i.e., $\numxi=1$).  The maximum polynomial degree used in the stochastic
space $S_{\numpsi}$ is $p$; thus, the dimension of $S_{\numpsi}$ is
$\numpsi = p +1$. In physical space, the SPDE is defined over a two-dimensional rectangular
bounded domain $D$, and it is discretized using the finite element
method with bilinear ($Q_1$) elements as implemented in the
Incompressible Flow and Iterative Solver Software (IFISS) package
\cite{ifiss}. Sixteen elements are employed in each dimension, leading to $n_x
= 225 = 15^2$ degrees of freedom excluding boundary nodes. All numerical experiments are performed on an Intel 3.1 GHz i7 CPU, 16 GB RAM using \textsc{Matlab} R2015a. 

\textbf{Measuring weighted $\ell^2$-norms}. For all LSPG methods, the weighted
$\ell^2$-norms can be measured by evaluating the expectations in the quadratic
form of the objective function shown in \eqref{eq:wlspg_quad_form}.
This requires evaluation of three expectations
\begin{equation}
\|M r(\bar{x}) \|_2^2 \defeq \bar{x}^T T_1 \bar{x} - 2 T_2^T \bar{x} + T_3,
\end{equation}
with 
\begin{align}
T_1 \defeq & E[ (\psi \psi^T \otimes A^TM^T M^TA)] \in \mathbb{R}^{n_x \numpsi \times n_x \numpsi},\label{eq:T1_exp}\\
T_2 \defeq & E[ \psi \otimes  A^T M^T  M \RHS{} ] \in \mathbb{R}^{n_x \numpsi},\label{eq:T2_exp}\\
T_3 \defeq & E[ \RHS{}^TM^T M\RHS{}] \in \mathbb{R}.\label{eq:T3_exp}
\end{align}
Note that $T_3$ does not depend on the stochastic-space dimension $\numpsi$.
These quantities can be evaluated by
numerical quadrature or analytically if closed-form expressions for those
expectations exist. Unless otherwise specified, we compute these quantities using the {\tt integral} function in \textsc{Matlab}, which performs
adaptive numerical quadrature based on the 15-point Gauss--Kronrod quadrature formula \cite{shampine2008vectorized}.

\textbf{Error measures}.
In the experiments, we assess the error in approximate solutions computed
using various spectral-projection techniques using 
four relative error measures (see Table \ref{tab:wlspg_method_table}):
\begin{align}\label{eq:error_measures}
\eta_r(x) \defeq\frac{\|\errorsol{x}\|_{A^TA}^2}{\|f\|_2^2},\quad
\eta_e(x) \defeq\frac{\|\errorsol{x}\|_2^2}{\|u\|_2^2},\quad
\eta_A(x) \defeq\frac{\|\errorsol{x}\|_A^2}{\|u\|_A^2},\quad
\eta_Q(x) \defeq\frac{\|F\errorsol{x}\|_2^2}{\|Fu\|_2^2}.
\end{align}

\subsection{Stochastic diffusion problems}\label{sec:diff}
Consider the steady-state stochastic diffusion equation with homogeneous boundary conditions,
\begin{equation}\label{eq:str_diff}
\left\{
\begin{array}{r l l}
-\nabla \cdot (a(x,\xi) \nabla u(x,\xi) ) &= f(x,\xi) &\text{ in } D \times \Gamma\\
u(x,\xi) &= 0 &\text{ on } \partial D \times \Gamma,
\end{array}
\right.
\end{equation}
where the diffusivity $a(x,\xi)$ is a random field and $D=[0, 1] \times [0,
1]$. The random field $a(x,\xi)$ is specified as an exponential of a truncated
Karhunen-Lo\`eve (KL) expansion \cite{loeve1978probability} with covariance kernel,
$C({x},\, {y}) \equiv \sigma^2 \exp \left (  - \frac{|x_1 - y_1|}{c} -
\frac{|x_2 - y_2|}{c} \right)$,
where $c$ is the correlation length, i.e.,
\begin{equation} \label{eq:def_rf}
a(x,\xi) \equiv \exp( \mu + \sigma a_1(x) \xi),
\end{equation} 
where $\{\mu, \sigma^2\}$ are the mean and variance of $a$ and $a_1(x)$ is the
first eigenfunction in the KL expansion. 
After applying the spatial
(finite-element) discretization, the problem can be reformulated as a parameterized linear system
of the form \eqref{eq:param_sys},
where 
$A(\xi)$ is a parameterized stiffness matrix obtained from the weak form of
the problem whose $(i,j)$-element is  $[A(\xi)]_{ij} = \int_D
\nabla a(x,\xi) \varphi_i(x) \cdot \varphi_j(x) dx$ (with $\{ \varphi_i \}$
standard finite element basis functions) and $\RHS{(\xi)}$ is a parameterized
right-hand side whose $i$th element is $[\RHS{(\xi)}]_i = \int_D f(x,\xi) \varphi_i(x)
dx$.
Note that $A(\xi)$ is symmetric positive definite for this problem; thus
\LSPGSG{} is a valid projection scheme (the Cholesky factorization
$A(\xi) = C(\xi)C(\xi)^T$ exists) and is equal to stochastic Galerkin
projection.
%
%

\textbf{Output quantities of interest}. We consider $n_o$ output
quantities of interest ($F(\xi)u(\xi) \in \RR{n_o}$) that are random linear functionals of the solution and $F(\xi)$ is of dimension $n_o \times n_x$  having the form:
\begin{itemize}
\item[(1)] $F_1(\xi) \defeq g(\xi)\times  G$ with $G\sim U (0,\,1)^{n_o \times n_x}$: each entry of $G$ is drawn uniformly from $[0, 1]$ and $g(\xi)$ is a scalar-valued function of $\xi$. The resulting output QoI, $F_1(\xi)u(\xi)$, is a vector-valued function of dimension $n_o$.
\item[(2)] $F_2(\xi)\defeq \RHS{(\xi)}^T \bar{M} $: $\bar{M}$ is a mass matrix defined via $[\bar{M}]_{ij} \equiv \int_D \varphi_i(x) \varphi_j(x)dx$. The output QoI is a scalar-valued function $F_2(\xi)u(\xi) = \RHS{(\xi)}^T \bar{M} u(\xi) $, which approximates a spatial average $\frac{1}{| D |}\int_D f(x,\xi) u(x,\xi) dx$. 
\end{itemize}

\begin{figure}[!htb]
	\centering
	\subfloat[Relative energy norm of solution error $\eta_A$] {
		\includegraphics[angle=0, scale=0.77]{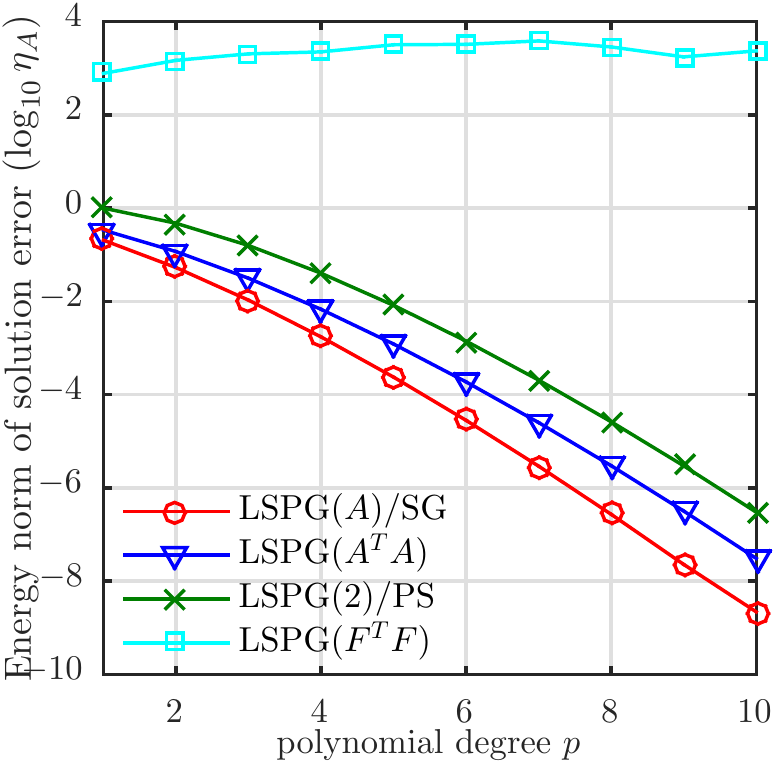}
		\label{fig:easy_pvsq_sub3}
	}\hspace{2mm}
	\subfloat[Relative $\ell^2$-norm of residual $\eta_r$] {
		\includegraphics[angle=0, scale=.77]{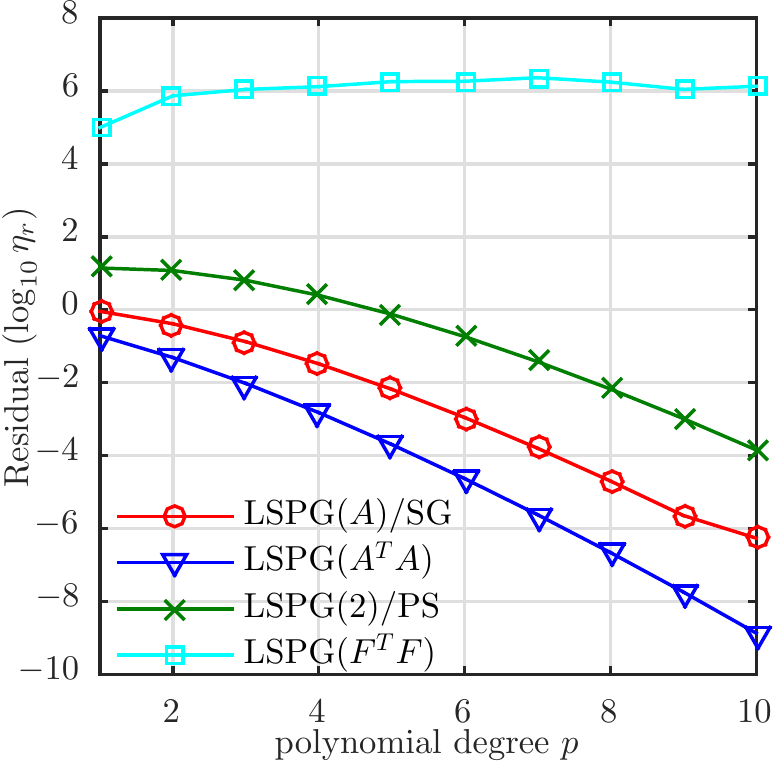}
		\label{fig:easy_pvsq_sub1}
	}\\
	\subfloat[Relative $\ell^2$-norm of solution error $\eta_e$] {
		\includegraphics[angle=0, scale=0.77]{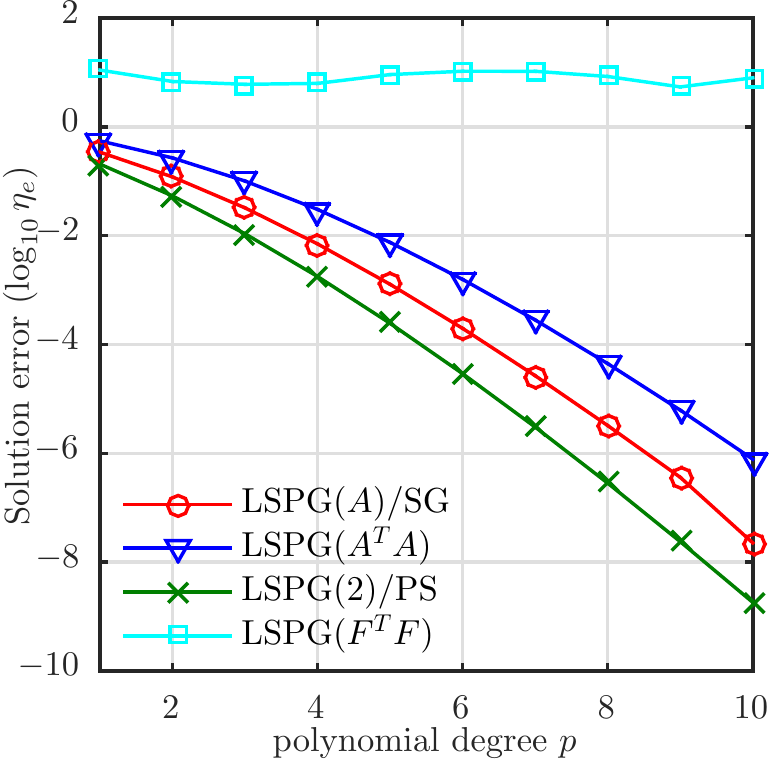}
	}
	\hspace{2mm}
	\subfloat[Relative $\ell^2$-norm of output QoI error $\eta_Q$ with $F=F_1$, $n_o=100$, and $g(\xi) =\xi$  ] {
		\includegraphics[angle=0, scale=0.77]{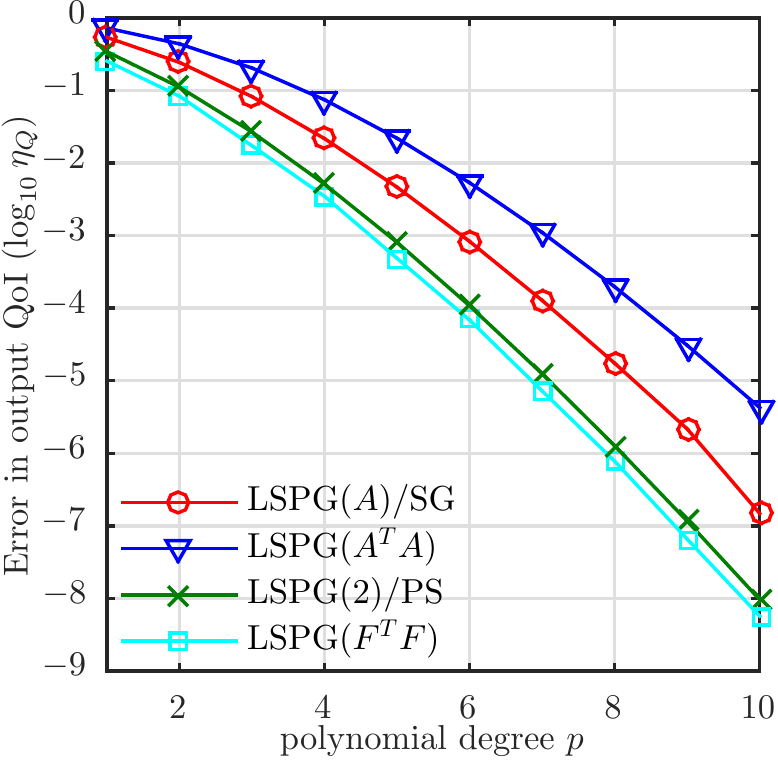}
	}

	\caption{Relative error measures versus polynomial
	degree for diffusion problem 1: lognormal random coefficient and
	deterministic forcing. Note that each LSPG method performs best in the error measure it
	minimizes.}
	\label{fig:easy_pvsq}
\end{figure} 


\subsubsection{Diffusion problem 1: Lognormal random coefficient and
deterministic forcing}\label{sec:easy_example}
In this example, we take $\xi$ in \eqref{eq:def_rf} to follow a standard normal distribution (i.e., $\rho(\xi) = \frac{1}{\sqrt{2\pi}}\exp\left(-\frac{\xi^2}{2}\right)$ and $\xi \in (-\infty, \infty)$) and $f(x,\xi) = 1$ is deterministic. Because $\xi$ is normally
distributed, normalized Hermite polynomials (orthogonal with respect to $\langle\cdot,\cdot\rangle_\rho$) are used as polynomial basis $\{\psi_i(\xi)\}_{i=1}^{\numpsi} $.

\begin{figure}[!tb]
	\hspace{0.1mm}
	\subfloat[Relative energy norm of solution error $\eta_A$] {
		\includegraphics[angle=0, scale=0.76]{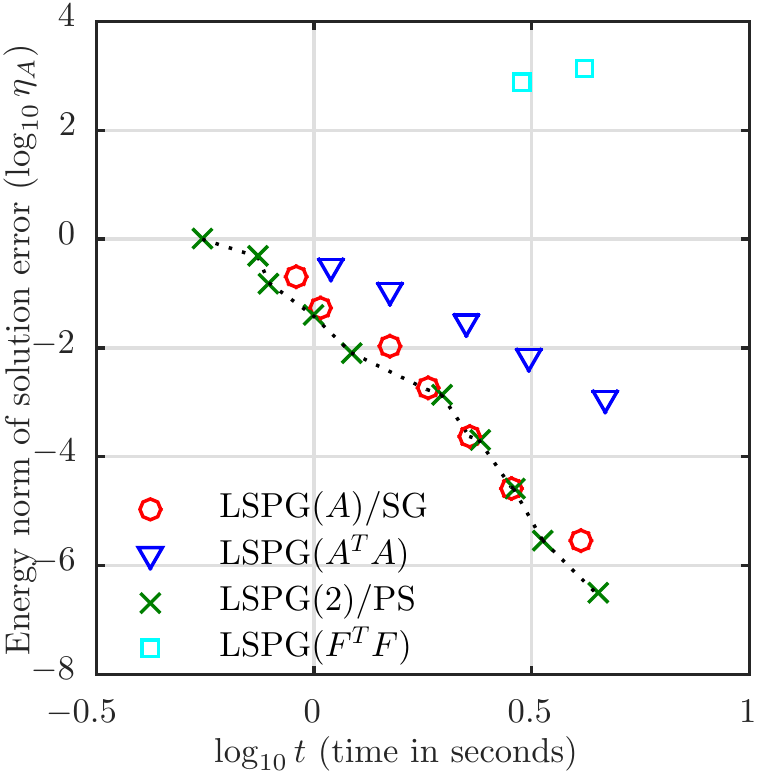}
	}
	\hspace{2.1mm}
	\subfloat[Relative $\ell^2$-norm of residual $\eta_r$] {
		\includegraphics[angle=0, scale=0.76]{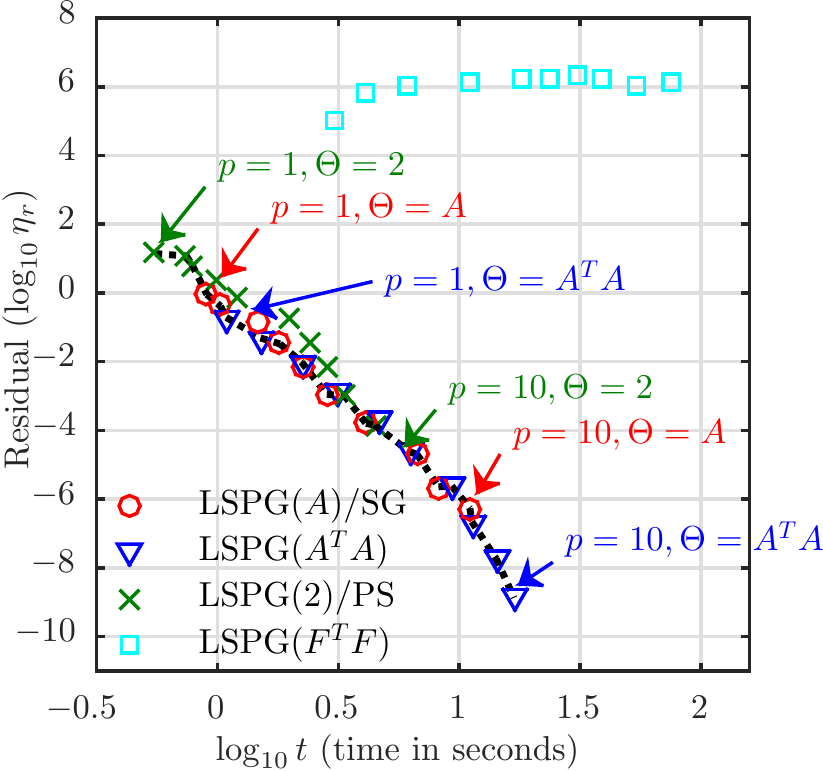}
	}\\
	\subfloat[Relative $\ell^2$-norm of solution error $\eta_e$] {
		\includegraphics[angle=0, scale=0.76]{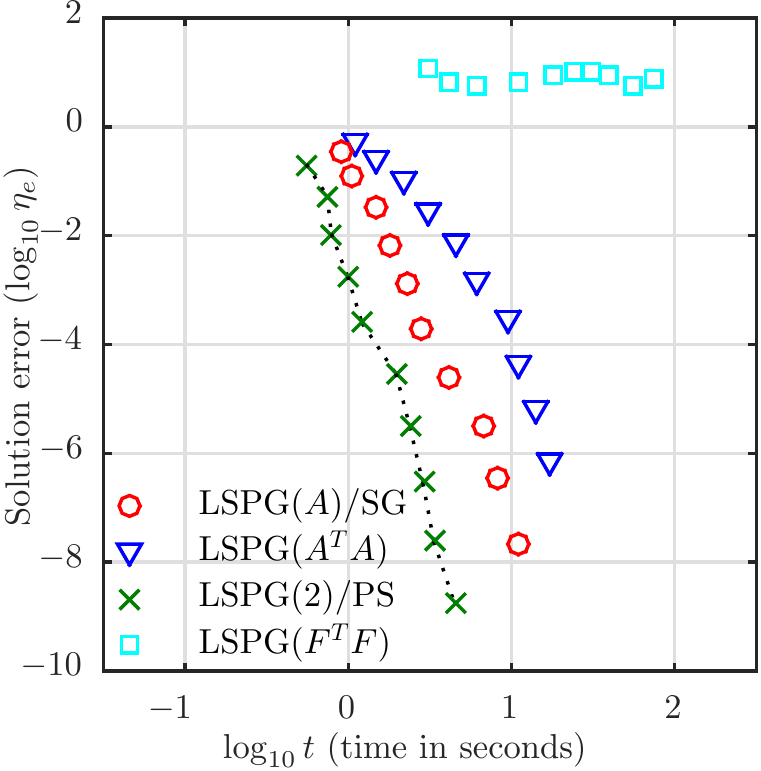}
	}
	\hspace{2mm}
	\subfloat[Relative $\ell^2$-norm of output QoI error $\eta_Q$ with $F=F_1$, $n_o=100$, and $g(\xi) =\xi$ ] {
		\includegraphics[angle=0, scale=0.76]{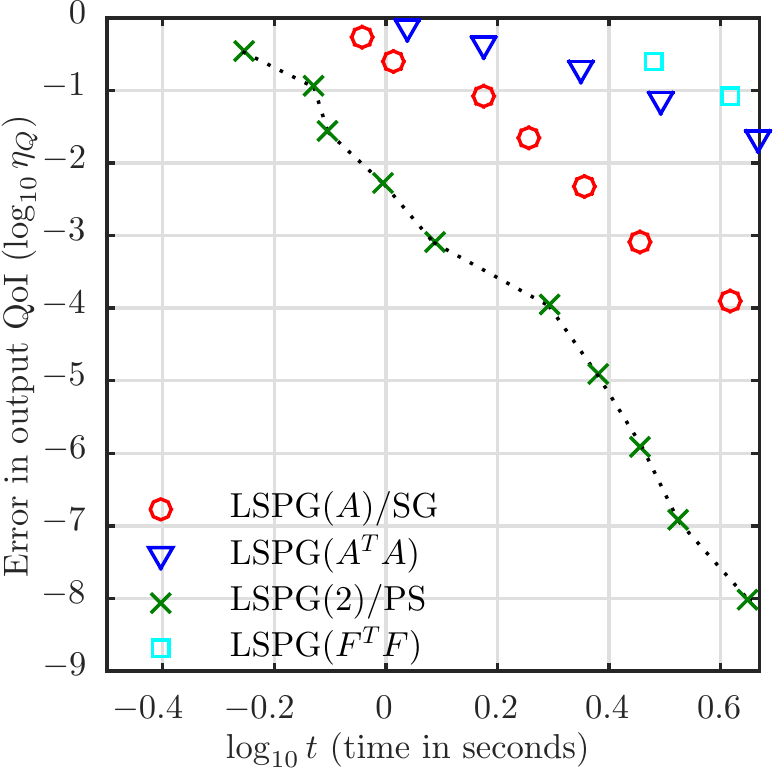}
	}\hspace{3.5mm}

	\caption{Pareto front of relative error measures versus wall time for varying polynomial degree $p$ for diffusion problem 1: lognormal random coefficient and deterministic forcing.}
	\label{fig:easy_quad}
\end{figure} 

Figure \ref{fig:easy_pvsq} reports the relative errors
\eqref{eq:error_measures} associated with solutions computed using four LSPG methods
(\LSPGSG{}, \LSPGI{}, \LSPGPS{}, and \LSPGQ{}) for varying
polynomial degree $p$. Here, we consider the random output QoI, i.e., $F=F_1$, $n_o=100$, and $g(\xi) = \xi$.
This result
shows that three methods (\LSPGSG{}, \LSPGI{}, and  \LSPGPS{})
monotonically converge in all four error measures, whereas \LSPGQ{} does
not. This is an artifact of rank deficiency in $F_1$, which leads to
$\sigma_\text{min}(F_1) = 0$; as a result, all stability constants $C$ for which
$\Theta = F^TF$ in Table \ref{tab:norm_equivalence} are unbounded, implying
lack of error control. 
Figure \ref{fig:easy_pvsq} also shows that each LSPG method
minimizes its targeted error measure for a given 
stochastic-subspace dimension (e.g., LSPG minimizes the $\ell^2$-norm of the residual); this is also evident from Table
\ref{tab:norm_equivalence}, as the stability constant realizes its minimum
value ($C=1$) for $\Theta = \Theta'$. Table \ref{tab:norm_equiv_for_diff1}
shows actual values of the stability constant of this problem and well
explains the behaviors of all LSPG methods. For example, the first column of
Table \ref{tab:norm_equiv_for_diff1} shows that the stability constant is
increasing in the order (\LSPGSG, \LSPGI, \LSPGPS, and \LSPGQ), which is
represented in Figure \ref{fig:easy_pvsq_sub3}. 

 \begin{table}[ht] 
 \begin{center}\footnotesize
\renewcommand{\arraystretch}{1.3}
\setlength{\tabcolsep}{4pt}
\caption{Stability constant $C$ of Diffusion problem 1}
 \renewcommand{\arraystretch}{1.2}
 \begin{tabular}{|l|r|r|r|r|} 
\hline
&$\Theta' = A$&$\Theta' = A^TA$&$\Theta' = 2$ & $\Theta' = F^TF$\\
\hline
$\Theta = A$ & 1	& 26.43 &	2.06 &	11644.22\\
$\Theta = A^TA$ & 2.06 &	1	& 4.25 & 	24013.48 \\
$\Theta = 1$ & 26.43 &	698.53 &	1 &	5646.32\\
$\Theta = F^TF$ & $\infty$ & $\infty$ & $\infty$ & 1\\
\hline 
\end{tabular} 
\label{tab:norm_equiv_for_diff1}
\end{center}
  \end{table}

The results in Figure \ref{fig:easy_pvsq} do not account for computational costs. This point is addressed in Figure \ref{fig:easy_quad}, which shows the relative errors as a function of CPU time.
As we would like to devise a method that minimizes both the error and computational time, we examine a Pareto front (black dotted line) in each error measure. For a fixed value of $p$,
\LSPGPS{} is the fastest method because it does not require solution of a coupled system of linear equations of dimension $n_x \numpsi$ which is required by the other three LSPG methods (\LSPGSG, \LSPGI, and \LSPGQ).
As a result, pseudo-spectral projection  (\LSPGPS) generally
yields the best overall performance in practice, even when it produces larger
errors than other methods for a fixed value of $p$.
Also, for a fixed value of $p$, \LSPGSG{} is faster than \LSPGI{} because the weighted stiffness matrix
$A(\xi)$ obtained from the finite element discretization is sparser 
than $A^T(\xi)A(\xi)$. That is, the number of nonzero entries to be
evaluated for \LSPGSG{} in numerical quadrature is smaller than the ones
for \LSPGI, and exploiting this sparsity structure in the numerical quadrature
	causes \LSPGSG{}  to be faster than \LSPGI. 

\subsubsection{Diffusion problem 2: Lognormal random coefficient and random
forcing}\label{sec:hard_example}
\begin{figure}[!htbp]
	\centering
	\includegraphics[angle=0, scale=.7]{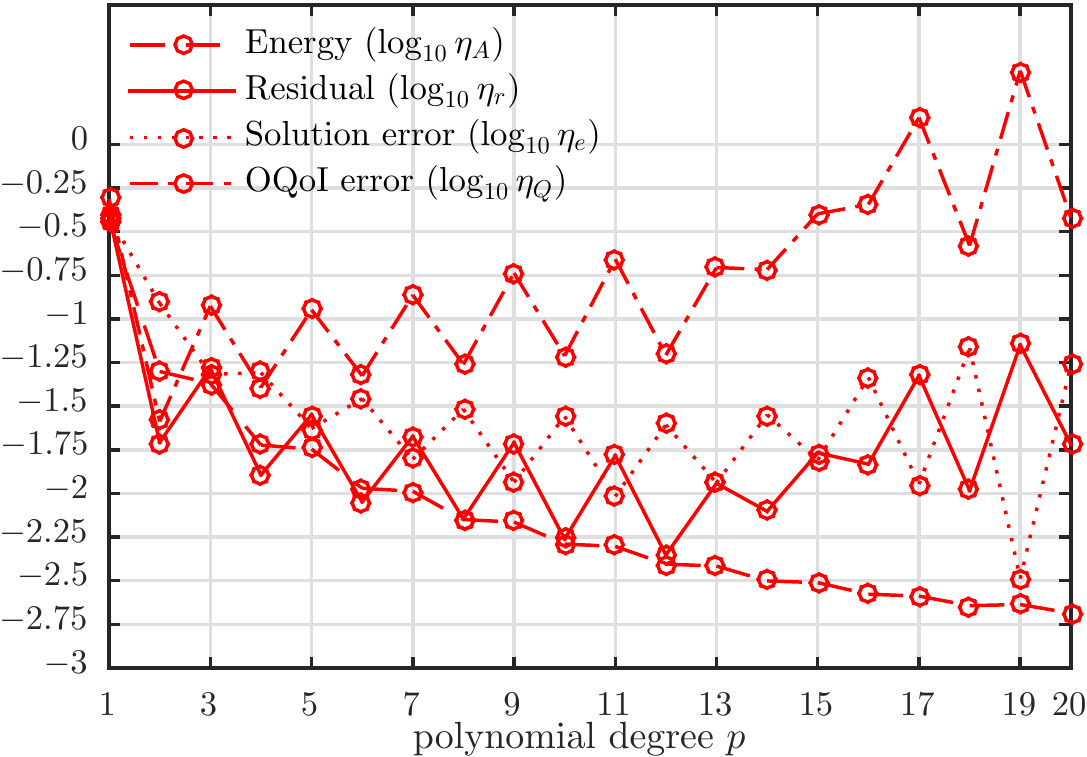}
	\caption{Relative errors versus polynomial degree for
	stochastic Galerkin (i.e., \LSPGSG) for diffusion problem 2: lognormal random coefficient and random forcing. Note that monotonic
	convergence is observed only in the minimized error measure $\eta_A$.}
	\label{fig:sg_fail}
\end{figure}

This example uses the same random field $a(x,\xi)$
\eqref{eq:def_rf}, but instead employs a random forcing term\footnote{In \cite{mugler2013convergence}, it was shown that stochastic Galerkin solutions of an analytic problem
$a(\xi)u(\xi)=f(\xi)$ with this type of forcing are divergent in the $\ell^2$-norm
of solution errors as $p$ increases.} $f(x,\xi) =
\exp(\xi)|\xi-1|$. Again, $\xi$ follows a standard normal distribution and
normalized  Hermite polynomials are used as polynomial basis. 
We consider the second output QoI, $F=F_2$.
As shown in Figure \ref{fig:sg_fail},
the stochastic Galerkin method fails to converge monotonically in three error
measures as the
stochastic polynomial basis is enriched. In fact, it exhibits monotonic
convergence only in the error measure it minimizes (for which monotonic
convergence is guaranteed). 

\begin{figure}[!h]
	\centering
	\subfloat[Relative energy norm of solution error $\eta_A$] {
		\includegraphics[angle=0, scale=0.76]{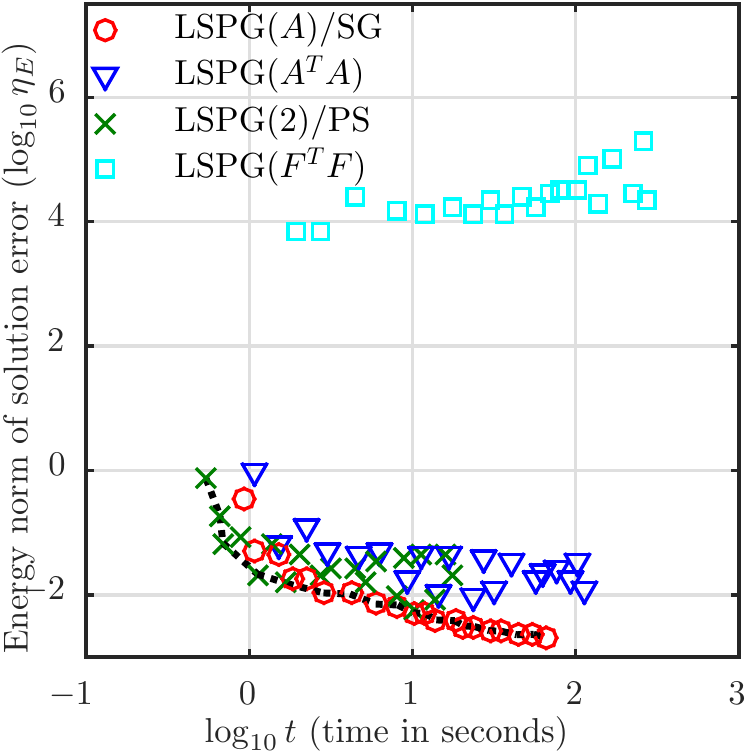}
	}\hspace{1mm}
	\subfloat[Relative $\ell^2$-norm of residual $\eta_r$] {
		\includegraphics[angle=0, scale=0.76]{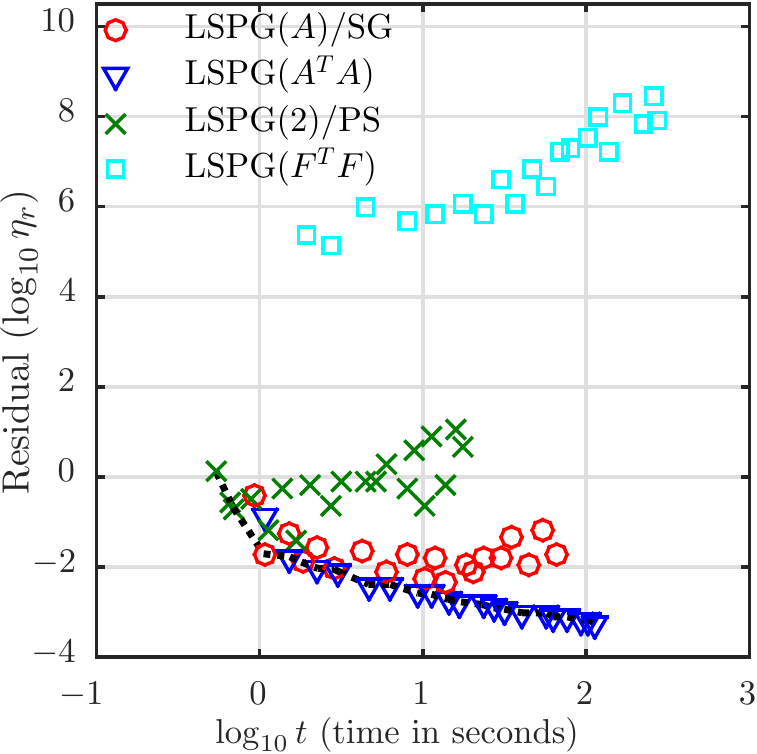}
	}\\
	\subfloat[Relative $\ell^2$-norm of solution error $\eta_e$] {
		\includegraphics[angle=0, scale=0.76]{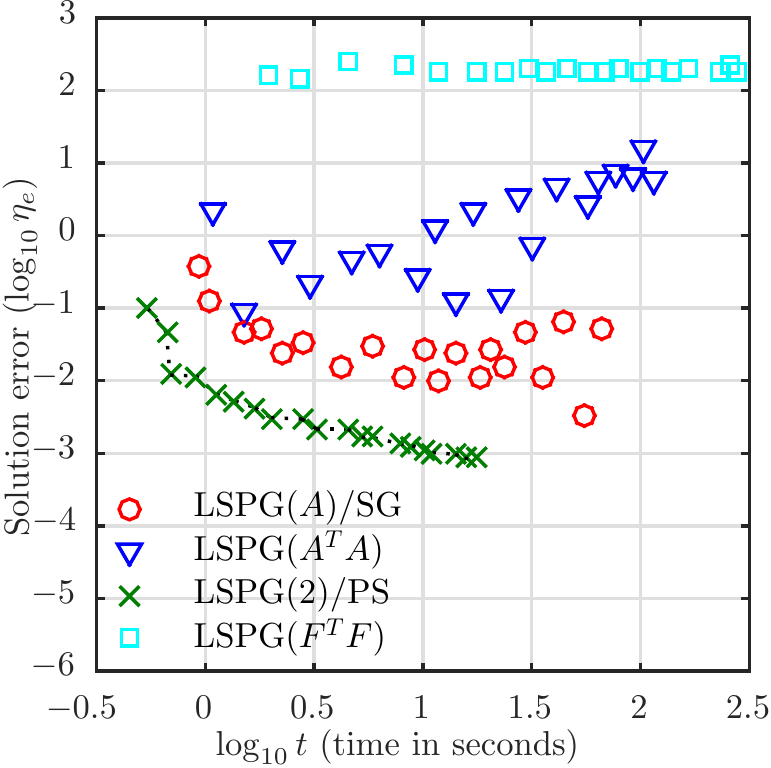}
	}\hspace{1mm}
	\subfloat[Relative $\ell^2$-norm of output QoI error $\eta_Q$ with $F=F_2$] {
		\includegraphics[angle=0, scale=0.76]{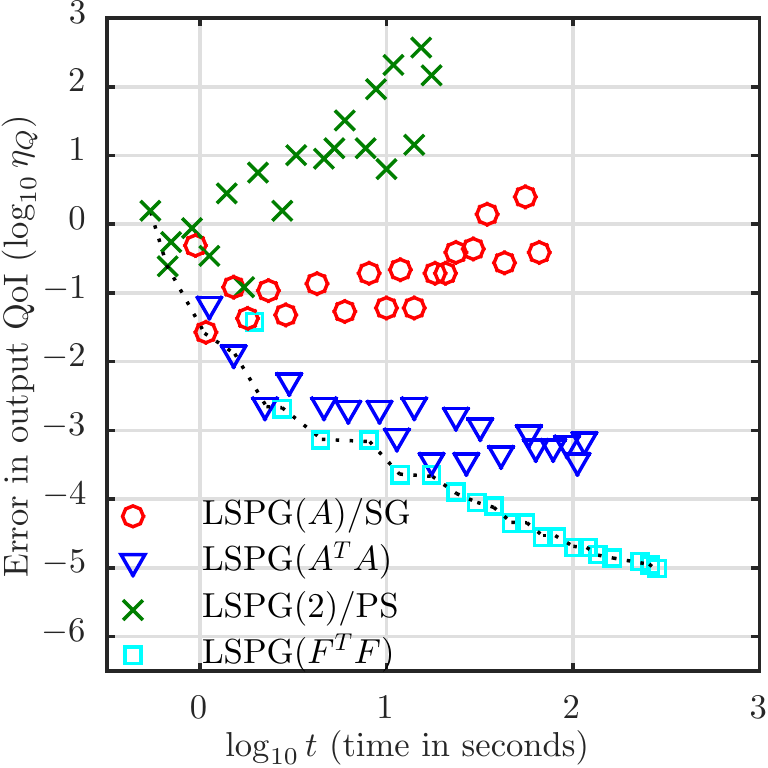}
	}

	\caption{Pareto front of relative error measures versus
	wall time for varying polynomial degree $p$ for diffusion problem 2: lognormal random coefficient and random forcing}
	\label{fig:hard_quad}
\end{figure} 

\begin{figure}[!t]
	\centering
	\subfloat[Relative energy norm of solution error $\eta_A$] {
		\includegraphics[angle=0, scale=0.76]{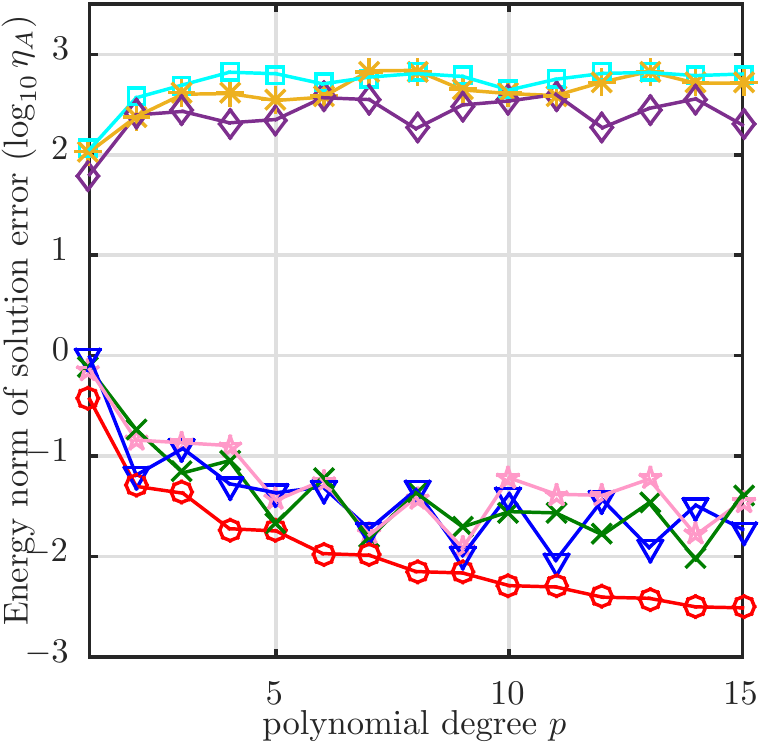}
	}\hspace{1mm}
	\subfloat[Relative $\ell^2$-norm of residual $\eta_r$] {
		\includegraphics[angle=0, scale=0.76]{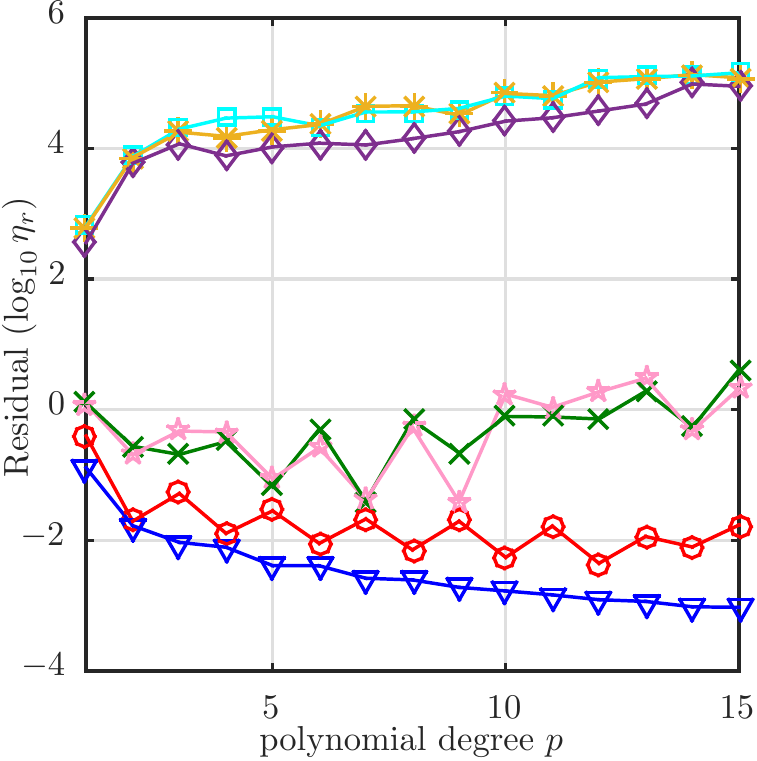}
	}\\
	\subfloat[Relative $\ell^2$-norm of solution error $\eta_e$] {
		\includegraphics[angle=0, scale=0.76]{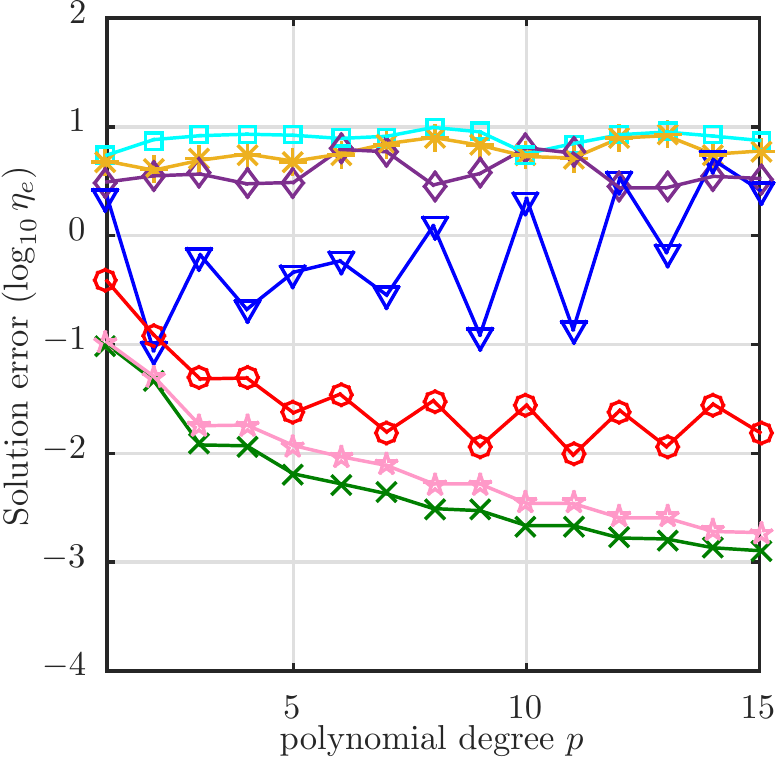}
	}\hspace{4.5mm}
	\subfloat[Legend for subplots (b)--(d)] {
		\includegraphics[angle=0, scale=.65]{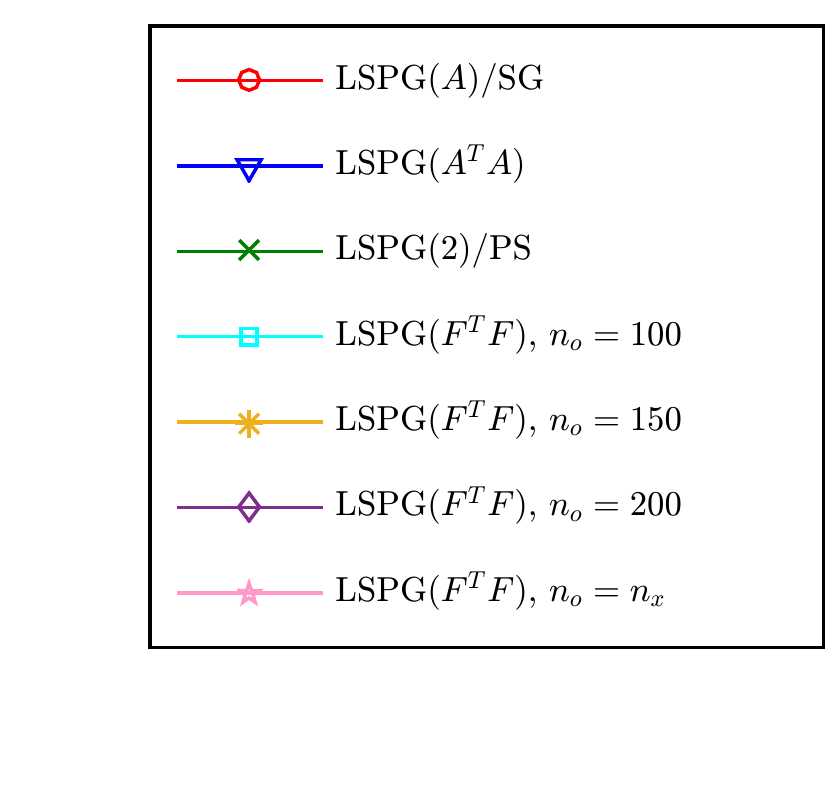}
	}
	\caption{Relative error measures versus polynomial degree for
	a varying dimension $n_o$ of the output matrix $F=F_1$ for diffusion problem 2: lognormal random coefficient and random forcing. Note	that \LSPGQ{} has controlled errors only when $n_o=n_x$, in which
	case $\sigma_\text{min}(F) > 0$.}
	\label{fig:hard_quad_qoi_three}
\end{figure}

\begin{figure}[!t]
	\centering
	\subfloat[
 Relative $\ell^2$-norm of output QoI error $\eta_Q$ with $F=F_1$, $n_o=100$, and
	$g(\xi) = g_1(\xi)=\sin(\xi)$] {
		\includegraphics[angle=0, scale=0.76]{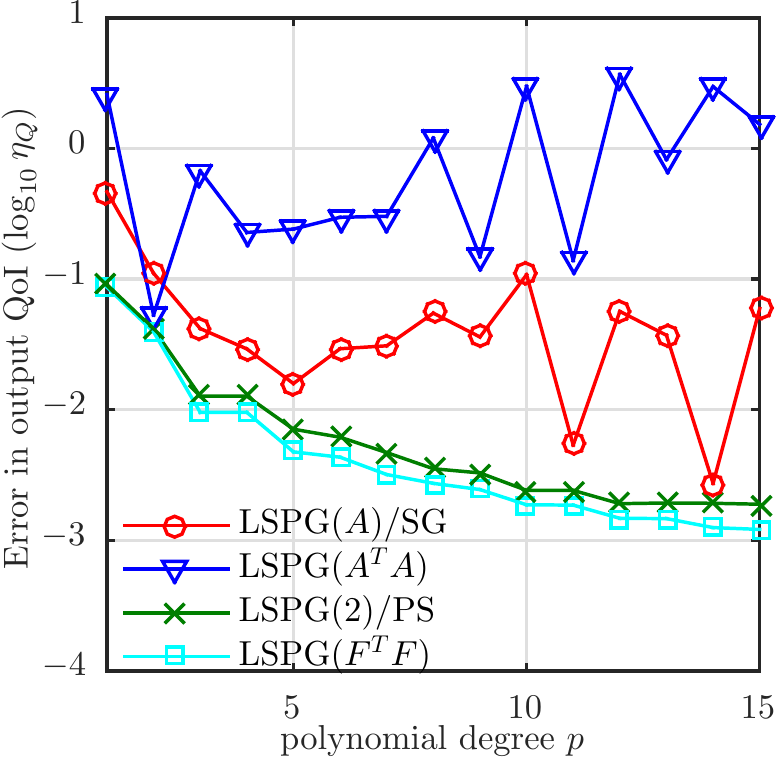}
	}
	\hspace{3mm}
	\subfloat[ Relative $\ell^2$-norm of output QoI error $\eta_Q$ with $F=F_1$, $n_o=100$, and
	$g(\xi) = g_2(\xi)=\xi$]
	 {
		\includegraphics[angle=0, scale=0.76]{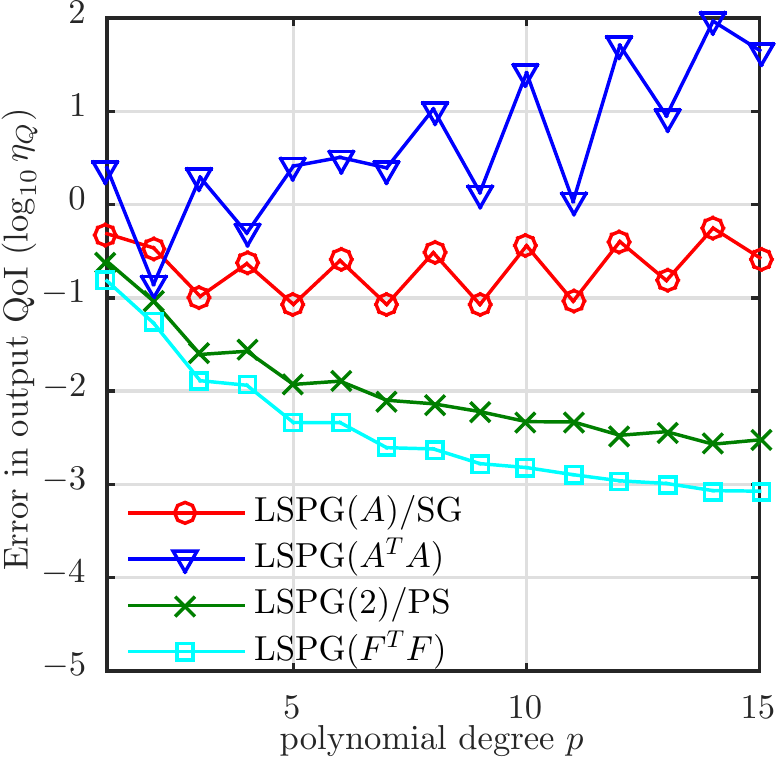}
	}
	\caption{Plots of the error norm of output QoI for diffusion problem 2: lognormal random coefficient and random forcing when a linear functional is (a) $F(\xi) \equiv \sin(\xi) \times [0,1]^{100\times n_x} $ and (b) $F(\xi) = \xi \times [0,1]^{100\times n_x}$ for varying $p$ and varying $n_o$.}
	\label{fig:hard_quad_qoi_qoi}
\end{figure}

Figure \ref{fig:hard_quad} shows that this trend applies to other methods as
well when effectiveness is viewed with respect to CPU time; each technique exhibits monotonic convergence in its tailored error
measure only. Moreover, the Pareto fronts (black dotted lines) in each
subgraph of Figure \ref{fig:hard_quad} shows that the LSPG method tailored for
a particular error measure is Pareto optimal in terms of minimizing the error
and computational wall time. In the next experiments, we examine goal-oriented \LSPGQ{} for varying number of output quantities of interest $n_o$ and its effect on the stability constant $C$. Figure \ref{fig:hard_quad_qoi_three} reports three error measures computed
using all four LSPG methods. For \LSPGQ, the first linear function
$F=F_1$ is applied with $g(\xi) = \sin(\xi)$ and a varying number of outputs
$n_o = \{100,
150, 200, 225\}$. When
$n_o=225$, \LSPGQ{} and \LSPGPS{} behave similarly in all three
weighted $\ell^2$-norms. This is because when $n_0=225=n_x$, then
$\sigma_\text{min}(F) > 0$, so the stability constants $C$ for $\Theta =
F^TF$ in Table \ref{tab:norm_equivalence} are bounded. Figure \ref{fig:hard_quad_qoi_qoi} reports relative errors in the quantity of
interest $\eta_Q$ associated with linear functionals $F=F_1$ for two different functions $g(\xi)$,
$g_1(\xi) = \sin(\xi)$ and $g_2(\xi) = \xi$. Note that \LSPGSG{} and
\LSPGI{} fail to converge, whereas \LSPGPS{} and \LSPGQ{} converge, which can be explained by the stability constant $C$ in Table \ref{tab:norm_equivalence} where $\LSV{A} = 26.43$ and $\SSV{A} = 0.48$ for the linear operator $A(\xi)$ of this problem. \LSPGQ{} converges monotonically and produces the smallest error (for a
fixed polynomial degree $p$) of all the methods as expected.

\subsubsection{Diffusion problem 3: Gamma random coefficient and random
forcing}\label{sec:gamma_example}
This section considers a stochastic diffusion problem parameterized by a random variable that has a Gamma distribution, where $a(x,\xi) \equiv \exp(1 + 0.25 a_1(x) \xi + 0.01 \sin(\xi))$
with density $\rho(\xi) \equiv \frac{\xi^{\alpha} \exp(-\xi)}{\bar{\Gamma}(\alpha+1)}$, $\bar{\Gamma}$ is the Gamma function, $\xi \in [0, \infty)$, and $\alpha=0.5$. Normalized Laguerre polynomials (which are orthogonal with respect to $\langle \cdot, \cdot \rangle_\rho$) are used as polynomial basis. We consider a random forcing
term $f(x,\xi) =
\log_{10}(\xi) |\xi - 1|$ and the second QoI $F(\xi) = F_2(\xi) = \RHS{(\xi)}^T 
\bar{M}$. Note that numerical quadrature is the only option for computing expectations arise in this problem.

\begin{figure}[!t]
	\centering
	\subfloat[Relative energy norm of solution error $\eta_A$] {
		\includegraphics[angle=0, scale=0.76]{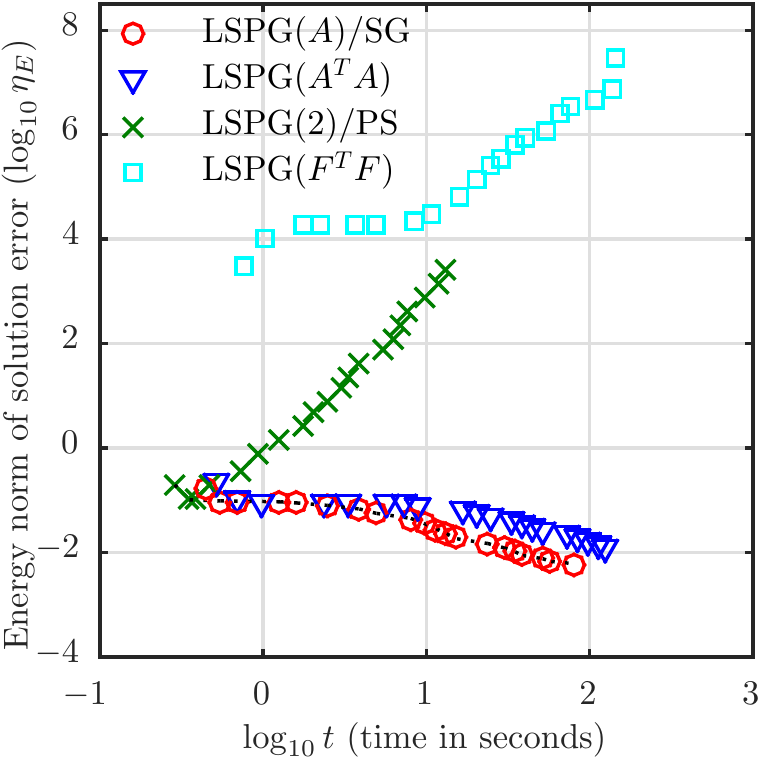}
	}\hspace{1mm}
	\subfloat[Relative $\ell^2$-norm of residual $\eta_r$] {
		\includegraphics[angle=0, scale=0.76]{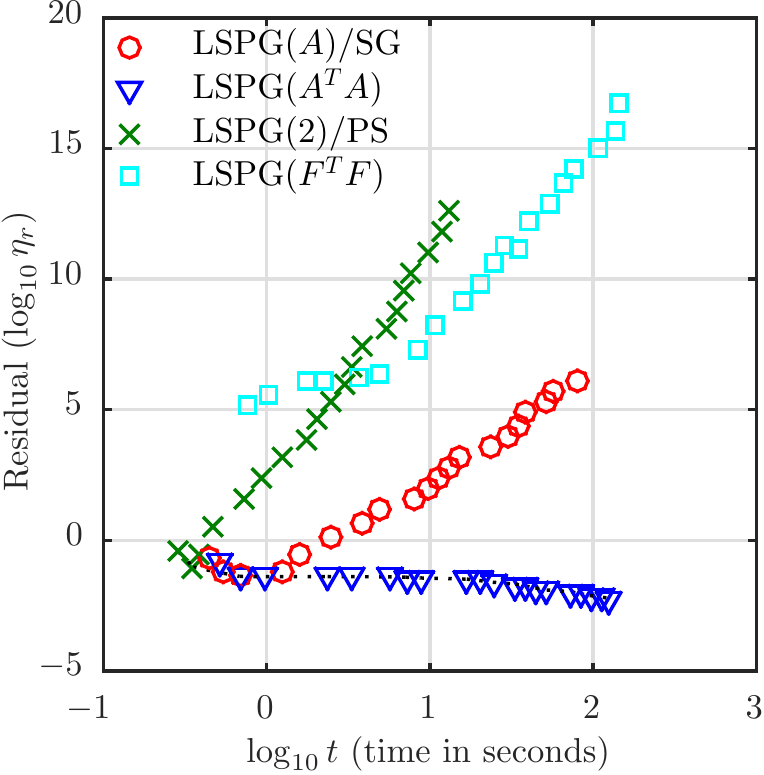}
	}\\
	\subfloat[Relative $\ell^2$-norm of solution error $\eta_e$] {
		\includegraphics[angle=0, scale=0.76]{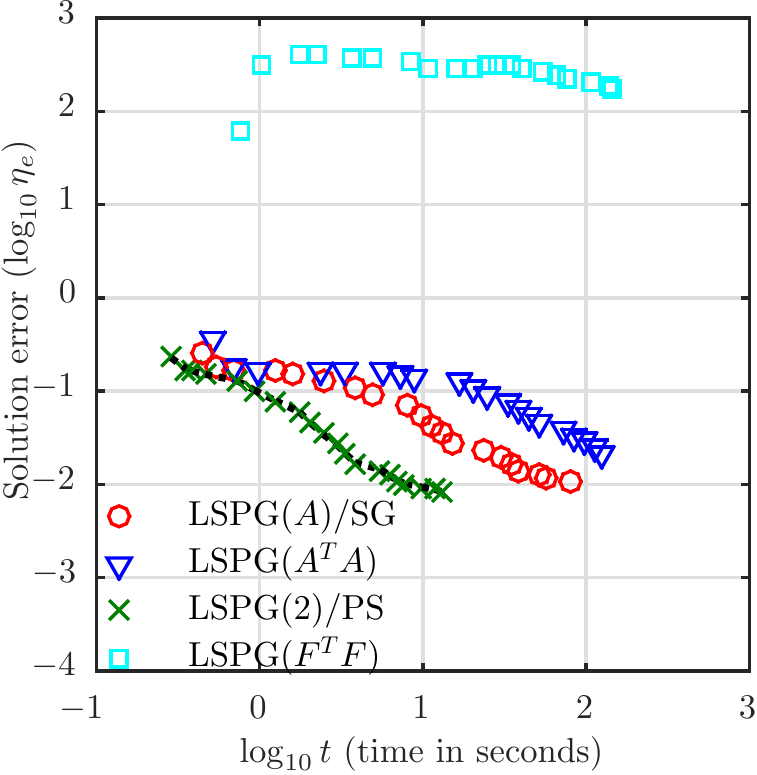}
	}\hspace{1mm}
	\subfloat[Relative $\ell^2$-norm of output QoI error $\eta_Q$ with $F=F_2$ ] {
		\includegraphics[angle=0, scale=0.76]{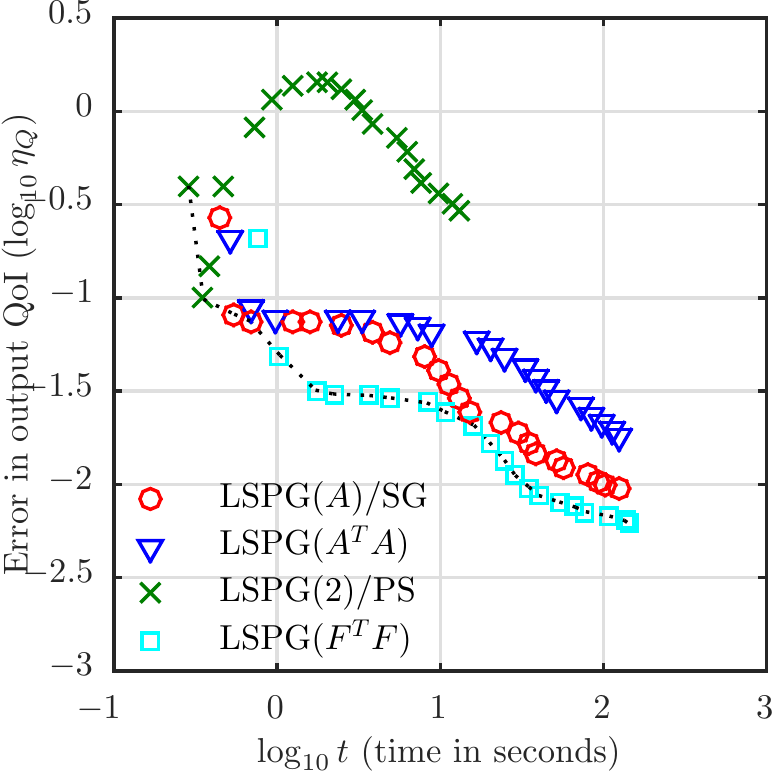}
	}

	\caption{Pareto front of relative error measures versus wall time for varying polynomial degree $p$ for diffusion problem 3: Gamma random coefficient and random
forcing. Note that each method is Pareto optimal in terms of minimizing its targeted error measure and computational wall time.}
	\label{fig:gamma_synthetic_quad}
\end{figure}

Figure \ref{fig:gamma_synthetic_quad} shows the
results of solving the problem with the four different LSPG methods. Again,
each version of LSPG monotonically decreases its corresponding target weighted
$\ell^2$-norm as the stochastic basis is enriched. Further, each LSPG method
is Pareto optimal in terms of minimizing its targeted error measure and the
computational wall time.

\subsection{Stochastic convection-diffusion problem: Lognormal random coefficient and
deterministic forcing}\label{sec:cd}
We now consider a non-self-adjoint example, the steady-state convection-diffusion equation 
\begin{equation}
\label{eq:cd_strong}
\left\{
\begin{array}{r l l}
-\epsilon \nabla \cdot (a({x},\,\xi) \nabla u({x},\,\xi) ) + \vec{w} \cdot \nabla u({x},\,\xi) &= f({x},\,\xi) &\text{ in } D \times \Gamma,\vspace{2mm}\\
u({x},\,\xi) &= g_D({x}) &\text{ on } \partial D \times \Gamma
\end{array}
\right.
\end{equation}
where $D=[-1,\,1]\times[-1,\,1]$ 
, $\epsilon$ is the viscosity parameter, and $u$ satisfies inhomogeneous Dirichlet boundary conditions 
\begin{equation}
g_D({x}) = 
\left\{
\begin{array}{c l}
g_D(x,\,1) = 0 &\text{ for } [-1, y] \cup [x,1] \cup [-1 \leq x \leq 0, -1],\\
g_D(1,\,y) = 1 &\text{ for } [1, y] \cup [0 \leq x \leq 1, -1].
\end{array}
\right.
\end{equation}
The inflow boundary consists of the bottom and the right portions of $\partial
D$, $[x,-1] \cup [1, y]$ \cite{elman2014finite}. We consider a zero forcing
term $f(x,\xi) =0$ and a constant convection velocity  $\vec{w} \equiv (-\sin \frac{\pi}{6}, \cos \frac{\pi}{6})$. We consider the convection-dominated case (i.e., $\epsilon = \frac{1}{200}$). 
\begin{figure}[!t]
	\centering
	\subfloat[Relative $\ell^2$-norm of residual $\eta_r$] {
		\includegraphics[angle=0, scale=0.76]{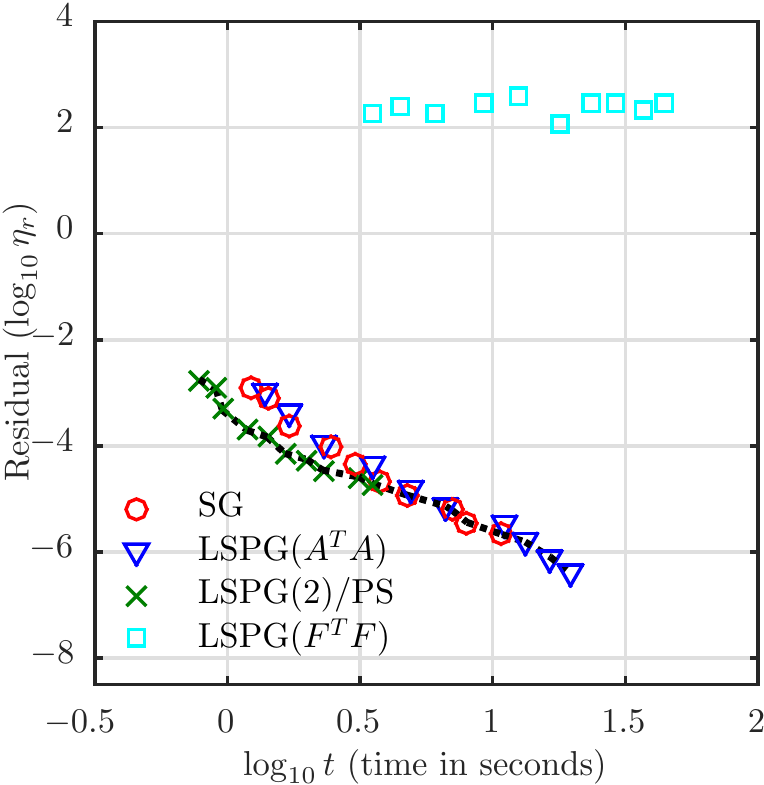}
	}\hspace{1mm}
	\subfloat[Relative $\ell^2$-norm of solution error $\eta_e$] {
		\includegraphics[angle=0, scale=0.76]{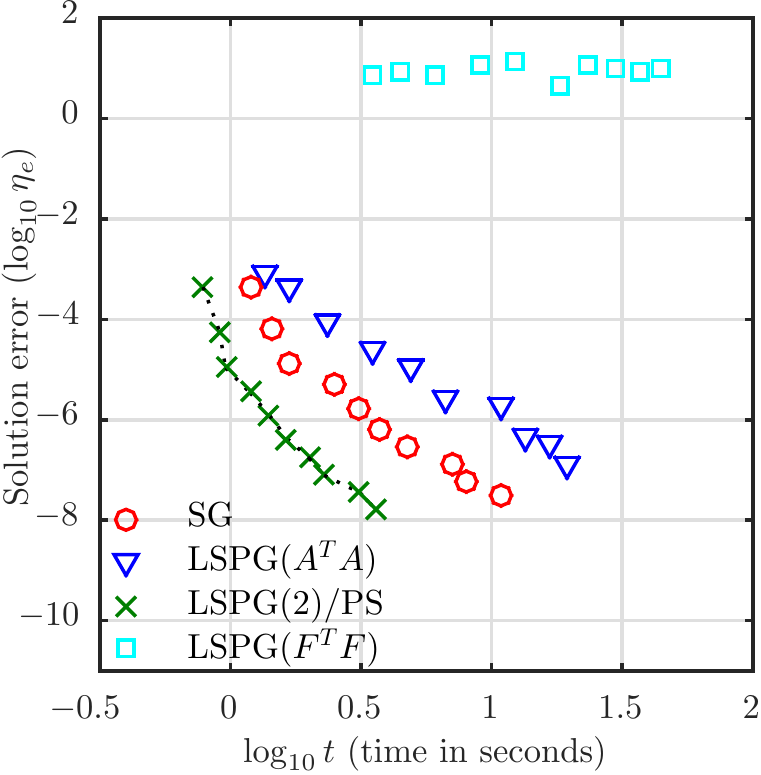}
	}\\
\subfloat[Relative $\ell^2$-norm of output QoI error $\eta_Q$ with $F=F_1$,
$n_o=100$, $g(\xi)=\exp(\xi) |\xi-1|$]
	 {
		\includegraphics[angle=0, scale=0.76]{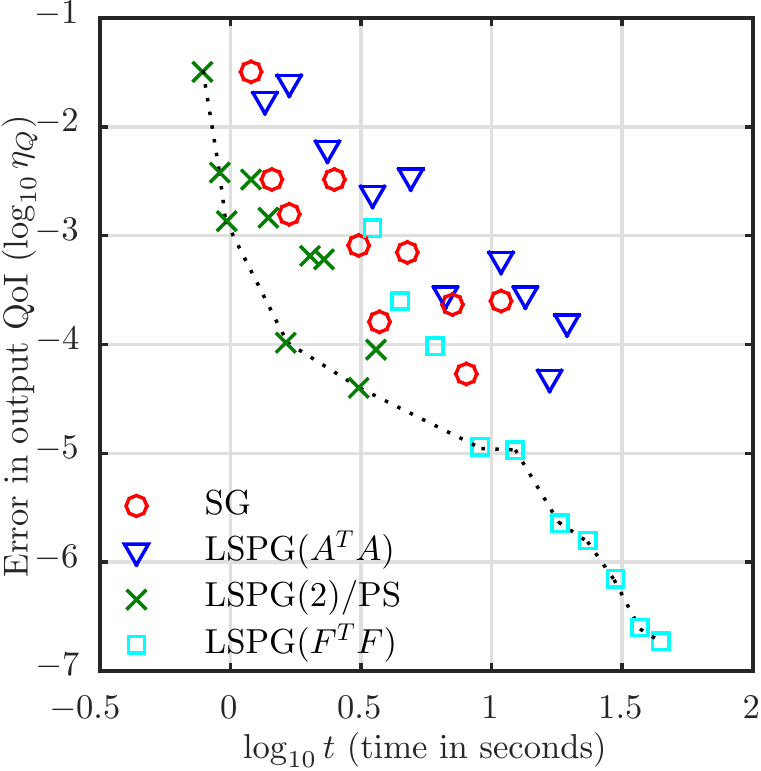}
	}

	\caption{Pareto front of relative error measures	versus wall time for varying polynomial degree $p$ for stochastic convection-diffusion problem: lognormal random coefficient and deterministic forcing term.}
	\label{fig:cd_easy}
\end{figure} 

For the spatial discretization, we essentially use the same finite element as above (bilinear $Q_1$ elements) applied to the weak formulation of \eqref{eq:cd_strong}. In addition, we use the streamline-diffusion method \cite{brooks1982streamline} to stabilize the discretization in elements with large mesh Peclet number. (See \cite{elman2014finite}, Ch. 8 for details.) Such spatial discretization leads to a parameterized linear system of the form \eqref{eq:param_sys} with
\begin{equation}\label{eq:param_sys_cd}
A(\xi) = \epsilon D(a(\xi); \xi) + C(\xi) + S(\xi),
\end{equation}
where $D(a(\xi);\xi)$, $C(\xi)$ and $S(\xi)$ are the diffusion term, the
convection term, and the streamline-diffusion term, respectively, and $[\RHS{(\xi)}]_i = \int_D f(x,\xi) \varphi_i (x) dx$. For this
numerical experiment, the number of  degrees of freedom in spatial domain is
$n_x=225$ (15 nodes in each spatial dimension) excluding boundary nodes. For \LSPGQ, the first linear function
$F=F_1$ is applied with $n_o=100$ outputs and $g(\xi) = \exp(\xi)|\xi-1|$.

Figure \ref{fig:cd_easy} shows the numerical results computed using the stochastic
Galerkin method and three LSPG methods (\LSPGI, \LSPGPS, \LSPGQ). Note that
the operator $A(\xi)$ is not symmetric positive-definite in this case; thus
LSPG($A$) is not a valid projection scheme (the Cholesky factorization $A(\xi)
= C(\xi)C(\xi)^T$ does not exist and the energy norm of the solution error
$\| \errorsol{x}\|_A^2$ cannot be defined) and stochastic
Galerkin does not minimize an any measure of the solution error. These results show that pseudo-spectral projection is Pareto optimal for achieving relatively larger error measures; this is because of its relatively low cost since, in contrast to the other methods,  it does not require the solution of a coupled linear system of dimension $n_x\numpsi$. In addition, the stochastic Galerkin projection is not Pareto optimal for any of the examples; this is caused by the lack of
optimality of stochastic Galekin in this case and highlights the significant benefit of
optimal spectral projection, which is offered by the stochastic LSPG method. In
addition, the residual $\eta_r$ and solution error $\eta_e$ incurred by
\LSPGQ{} are uncontrolled, because $n_o< n_x$ and thus $\sigma_\text{min}(F) = 0$.
Finally, note that each LSPG method is Pareto optimal for small errors in
its targeted error measure.

\subsection{Numerical experiment with analytic computations}
For the results presented above, expected values were computed using numerical
quadrature (using the \textsc{Matlab} function {\tt integral}). This is a
practical and general approach for numerically computing the required
integrals of \eqref{eq:T1_exp}--\eqref{eq:T3_exp}, and is the only option when
analytic computations are not available (as in Section
\ref{sec:gamma_example}). In this section, we briefly discuss how the costs
change if analytic methods based on closed-form integration exist and are used
for these integrals. Note that in general, however, analytic computation are
	unavailable, for example, if the random variables have a finite support (e.g., truncated Gaussian random variables as shown in \cite{ullmann2012efficient}). 


\textbf{Computing $T_1$}. Analytic computation of $T_1$ is possible if either $E[A^TMMA \psi_l]$ or $E[MA \psi_l]$ can be evaluated analytically. For \LSPGSG{} and \LSPGI, if $E[A \psi_l]$ can be evaluated so that the following gPC expansion can be obtained analytically
\begin{equation}\label{eq:pce_for_rf}
A(\xi) =  \sum_{l=1}^{\infty} A_l \psi_l (\xi), \quad A_l \equiv E[ A \psi_l],
\end{equation}
where $A_l \in \mathbb{R}^{n_x \times n_x}$, then $T_1$ can be computed analytically.  Replacing $A(\xi)$ with the series of \eqref{eq:pce_for_rf} for \LSPGSG{} ($M(\xi)=C^{-1}(\xi)$)  and \LSPGI{} ($M(\xi)=I_{n_x}$) yields 
\begin{align}\label{eq:T1_sg}
T_1\SGsup = \sum_{l=1}^{n_a}  E \left[ \psi\psi^T \otimes \left( A_l \psi_l \right)\right]  &= \sum_{l=1}^{n_a}  E [ \psi\psi^T \psi_l \otimes  A_l ],
\end{align}
and
\begin{align}\label{eq:T1_lspg}
T_1\LSPGsup = E[ \psi \psi^T \otimes \sum_{k=1}^{n_a} \sum_{l=1}^{n_a}\left(
A_k \psi_k \right)^T \left( A_l \psi_l \right)] = \sum_{k=1}^{n_a}
\sum_{l=1}^{n_a} E[ \psi \psi^T \psi_k \psi_l \otimes  A_k^T A_l ],
\end{align}
where the expectations of triple or quadruple products of the polynomial basis (i.e., $E[\psi_i \psi_j \psi_k]$ and $E[\psi_i \psi_j \psi_k \psi_l]$) can be computed analytically. 
For \LSPGPS, an analytic computation of $T_1$ is straightfoward because $M(\xi) A(\xi) = I_{n_x}$ and, thus,
\begin{align}
T_1\PSsup = E[ \psi\psi^T \otimes I_{n_x} ]  = I_{n_x \numpsi}.
\end{align}
Similarly, analytic computation of $T_1$ is possible for \LSPGQ if there
exists a closed formulation for $E[F \psi_l]$ or $E[F^TF \psi_l]$, which is again in general not available. 

\textbf{Computing $T_2$}. Analytic computation of $T_2$ can be performed in a similar way. If the random function $\RHS{(\xi)}$ can be represented using a gPC expansion, 
\begin{equation}\label{eq:f_pce}
\RHS(\xi) = \sum_{l=1}^{n_b} \RHS{}_l \psi_l(\xi), \quad \RHS{}_l \equiv E[ \RHS{} \psi_l ],
\end{equation}
 then, for \LSPGSG{} and \LSPGI, $T_2$ can be evaluated analytically by computing
 expectations of bi or triple products of the polynomial bases (i.e.,
 $E[\psi_i \psi_j ]$ and $E[\psi_i \psi_j \psi_k]$). For \LSPGPS{} and \LSPGQ,
 however, an analytic computation of $T_2$ is typically unavailable because a closed-form expression for $A^{-1}(\xi)$ does not exist. 

\begin{figure}[!h]
	\centering
	\subfloat[Relative energy norm of solution error $\eta_A$] {
		\includegraphics[angle=0, scale=0.76]{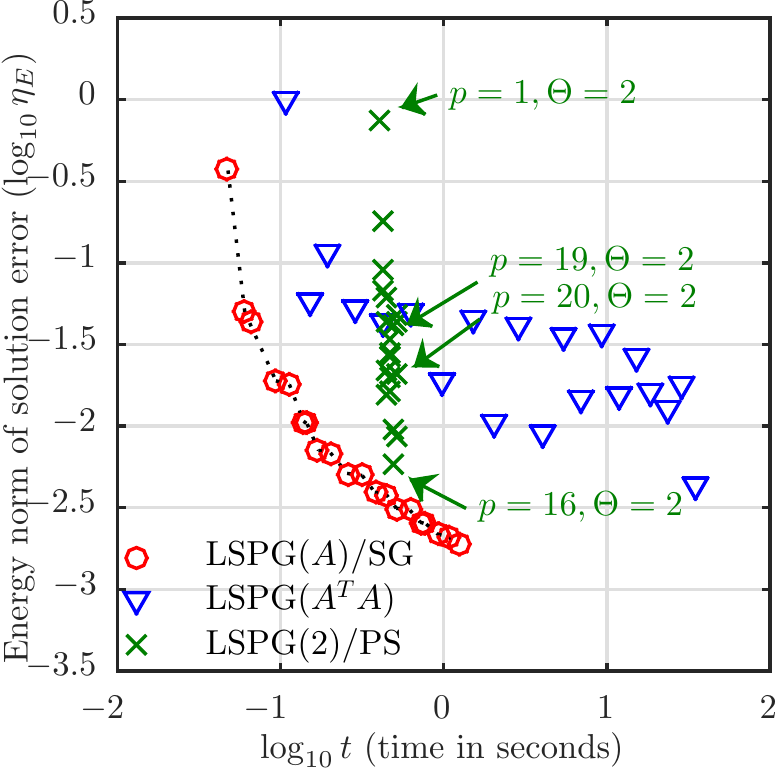}
	}\hspace{4mm}
	\subfloat[Relative $\ell^2$-norm of residual $\eta_r$]{
		\includegraphics[angle=0, scale=0.76]{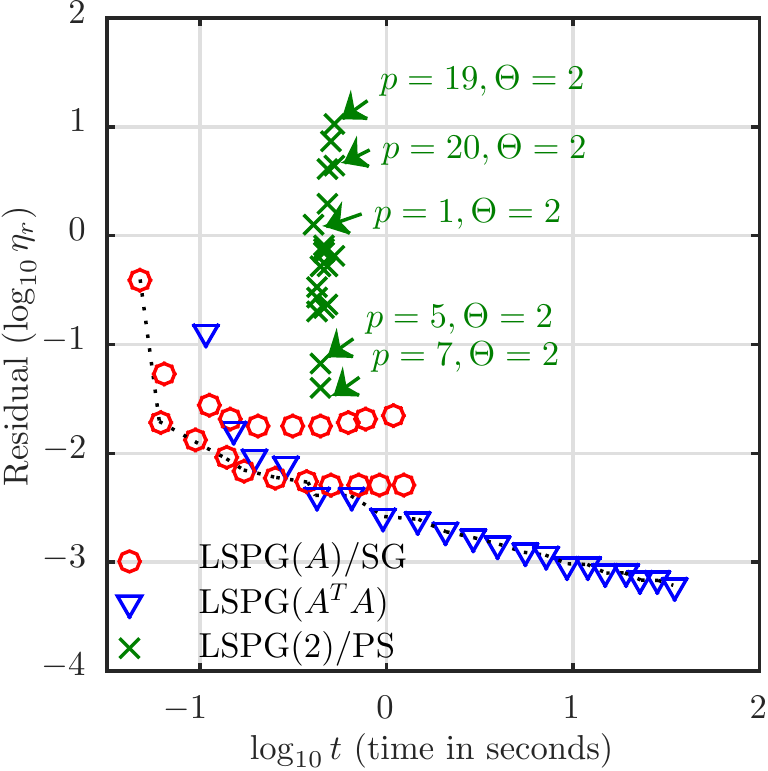}
	}\\
	\subfloat[Relative $\ell^2$-norm of solution error $\eta_e$] {
		\includegraphics[angle=0, scale=0.76]{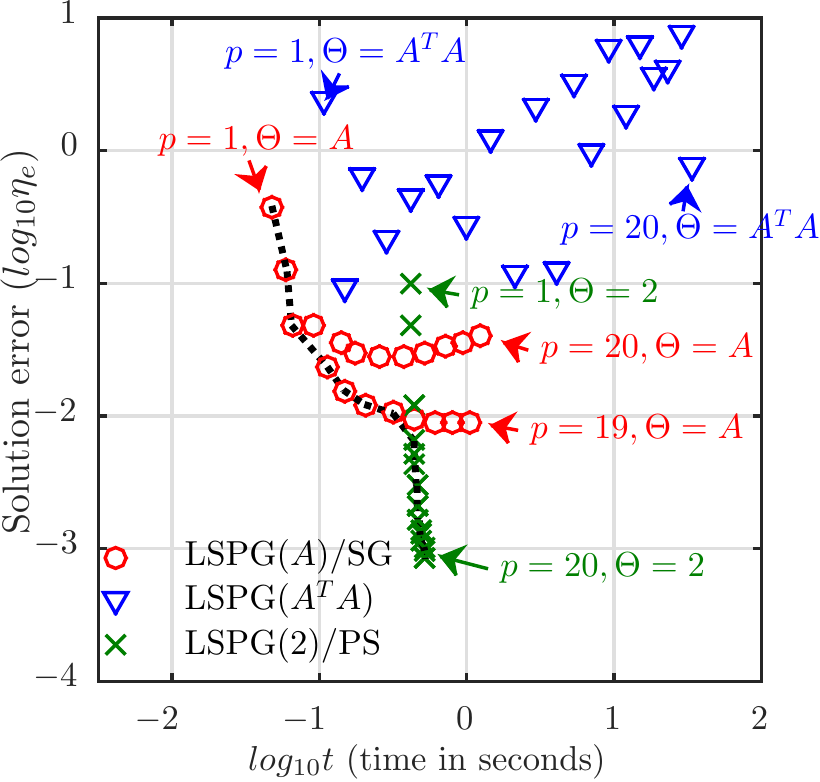}
	}\hspace{1mm}
	\subfloat[Relative $\ell^2$-norm of output QoI error $\eta_Q$ with $F = F_2$ ] {
		\includegraphics[angle=0, scale=0.76]{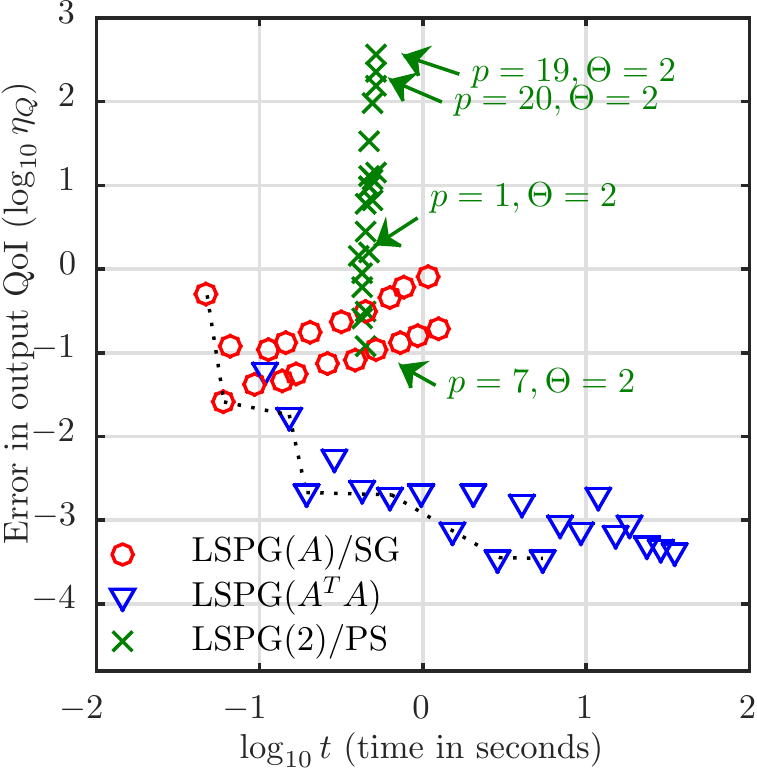}
	}

	\caption{Pareto front of relative error measures	versus wall time for varying polynomial degree $p$ for diffusion problem 2: Lognormal random coefficient and random
forcing. Analytic computations are used as much as possible to evaluate expectations.}
	\label{fig:hard_analytic}
\end{figure} 


We examine the impact of these observations on the cost of solution of the problem studied in Section
\ref{sec:hard_example}, a the steady-state stochastic diffusion equation \eqref{eq:str_diff} with lognormal random field $a(x,\xi)$ as in \eqref{eq:def_rf}, and random forcing $f(x,\xi) = \exp(\xi)|\xi-1|$.

Figure \ref{fig:hard_analytic} reports results for this problem for analytic computation of expectations. For \LSPGSG{}, analytic computation of the expectations $\{T_i\}_{i=1}^3$ requires fewer terms than for \LSPGI. In fact, comparing \eqref{eq:T1_sg} and \eqref{eq:T1_lspg} shows that computing $T_1\LSPGsup$ requires computing and assembling $n_a^2$ terms, whereas computing $T_1\SGsup$ involves only $n_a$ terms. Additionally the quantities $\{ A_k^T A_l \}_{k,l = 1}^{n_a}$ appearing in the terms of $T_1\LSPGsup$ in \eqref{eq:T1_lspg} are typically denser than the counterparts $\{A_k\}_{k=1}^{n_a}$  appearing in \eqref{eq:T1_sg}, as the sparsity pattern of $\{A_k\}_{k=1}^{n_a}$ is identical to that of the finite element stiffness matrices. As a result, \LSPGSG{} is Pareto optimal for small computational wall times when any error metric is considered.  When the polyomial degree $p$ is small, \LSPGSG{} is computationally faster than \LSPGPS, as \LSPGPS{} requires the solution
of $A(\xi^{(k)})u(\xi^{(k)}) = f(\xi^{(k)})$ at each quadrature point and cannot exploit analytic computation. As the stochastic basis is enriched, however, each tailored LSPG method outperforms other LSPG methods in minimizing its corresponding target error measure.

\section{Conclusion} \label{sec:conclusion}
In this work, we have proposed a general framework for optimal spectral projection
wherein the solution error can be minimized in weighted $\ell^2$-norms of
interest. In particular, we propose two new methods that minimize the
$\ell^2$-norm of the residual and the $\ell^2$-norm of the error in an output
quantity of interest. Further, we showed that when the linear operator is symmetric positive
definite, stochastic Galerkin is a
particular instance of the proposed methodology for a specific choice of
weighted $\ell^2$-norm. Similarly, pseudo-spectral projection is a particular case of
the method for a specific choice of weighted $\ell^2$-norm.

Key results from the numerical experiments include:
\begin{itemize} 
\item For a fixed stochastic subspace, each LSPG method minimizes its targeted
error measure (Figure \ref{fig:easy_pvsq}).
\item For a fixed computational cost, each LSPG method often minimizes its
targeted error measure (Figures \ref{fig:hard_quad},
\ref{fig:gamma_synthetic_quad}). However, this does
not always hold, especially for smaller computational costs (and smaller
stochastic-subspace dimensions) when larger errors
are acceptable. In
particular
pseudo-spectral projection (\LSPGPS) is often significantly less
expensive than other methods for a fixed
stochastic subspace, as it does not require solving a coupled linear system of
dimension $n_x\numpsi$ (Figures \ref{fig:easy_quad}, \ref{fig:cd_easy}). Alternatively, when analytic computations are possible,
stochastic Galerkin (\LSPGSG)) may be significantly less expensive than other methods for a fixed
stochastic subspace (Figure \ref{fig:hard_analytic}).
\item Goal-oriented \LSPGQ{} can have uncontrolled errors in error
measures that deviate from the output-oriented error measure $\eta_Q$ when the
linear operator $F$ has more columns $n_x$ than rows $n_o$ (Figure
\ref{fig:hard_quad_qoi_three}). This is because the minimum singular value
is zero in this case
(i.e., $\sigma_\text{min}(F)=0)$), which leads to unbounded stability
constants in other error measures (Table \ref{tab:norm_equivalence}).
\item Stochastic Galerkin often leads to divergence in different error
measures (Figure \ref{fig:sg_fail}). In this case, applying LSPG with the appropriate
targeted error measure can significantly improve accuracy (Figure
\ref{fig:hard_quad}).
\end{itemize}
Future work includes developing efficient sparse solvers for the stochastic LSPG methods and extending the methods to parameterized nonlinear systems.

\bibliography{references}
\bibliographystyle{siamplain}
\end{document}